\documentclass[11pt]{amsart}
\usepackage{amsmath,amsfonts,amsthm,amscd,amssymb,mathrsfs,amssymb,bm}
\usepackage{mathrsfs}
\usepackage[dvips]{graphicx}
\usepackage{wrapfig}
\usepackage[all]{xy}
\newtheorem{theorem}{Theorem}

\newtheorem{proposition}[theorem]{Proposition}

\newtheorem{remark}[theorem]{Remark}
\newcommand{\CP}{\mathbb{CP}}

\newcommand{\ol}{\overline}
\newcommand{\lra}{\longrightarrow}

\newcommand{\set}{\,|\,}
\newcommand{\proofend}{\hfill$\square$}

\newcommand{\inv}{^{-1}}
\newcommand{\Bs}{{\rm{Bs}}}

\newcommand{\Aut}{{\rm{Aut}}}

\newcommand{\vsp}{\vspace{3mm}}
\newcommand{\hsp}{\hspace{2mm}}
\numberwithin{equation}{section}
\numberwithin{theorem}{section}
\begin{document}
\bibliographystyle{alpha} 
\title{Geometry of generic Moishezon twistor spaces
on $4\mathbb{CP}^2$}
\author{Nobuhiro Honda}
\address{Mathematical Institute, Tohoku University,
Sendai, Miyagi, Japan}
\email{honda@math.tohoku.ac.jp}
\begin{abstract}
In this paper we investigate a family of Moishezon twistor spaces on the connected sum of 4 complex projective planes, which can be regarded as a direct generalization of the  twistor spaces on $3\mathbb{CP}^2$ of double solid type studied by Poon and Kreussler-Kurke.
These twistor spaces have a natural structure 
of double covering over a scroll of 2-planes
over a conic.
We determine the defining equations of the branch divisors in an explicit form,
which are very similar to the case of $3\mathbb{CP}^2$.
Using these explicit description we compute the 
dimension of the moduli spaces of these twistor spaces.
Also we observe that  
similarly to the case of $3\mathbb{CP}^2$, these twistor spaces can also be considered as generic Moishezon twistor spaces on $4\mathbb{CP}^2$.
We obtain these results by analyzing the anticanonical map of the twistor spaces in detail, which enables us to give an explicit construction of the twistor spaces, up to small resolutions.
\end{abstract}
\maketitle
\setcounter{tocdepth}{1}
\vspace{-5mm}


\section{Introduction}
In their papers, Kreussler-Kurke \cite{KK92} and Poon \cite{P92} investigated algebraic structure of generic twistor spaces on $3\mathbb{CP}^3$, the connected sum of 3 copies of complex projective planes.
They showed that if the half-anticanonical system of a twistor space of $3\mathbb{CP}^2$ is base point free, then the morphism associated to the system becomes a generically 2 to 1 covering map whose branch divisor is a quartic surface.
Further, they determined defining equation of the quartic surface; for the most generic twistor spaces,
 with respect to  homogeneous coordinates on $\mathbb{CP}^3$, the equation is of the form
\begin{align}\label{KKP}
z_0z_1z_2z_3=Q(z_0,z_1,z_2,z_3)^2
\end{align}
where $Q$ is a (homogeneous) quadratic polynomial with real coefficients. 
From the equation, the intersection of the quadratic surface $Q=0$ and any of the 4 plane $z_i=0$ is a double conic, and the intersection points of these 4 conics (consisting of 12 points) are ordinary double points of the quartic surface \eqref{KKP}.
For a generic quadratic polynomial $Q$,  these are all singularities of the surface \eqref{KKP}, but
they showed that when \eqref{KKP} is actually the branch divisor of the twistor spaces,  the surface has one more node, which is necessarily real, so that the branch surface has 13 ordinary nodes in total.
%

In this paper, we shall find Moishezon twistor spaces on $4\mathbb{CP}^2$ which can be regarded as a direct generalization
 of these twistor spaces on $3\mathbb{CP}^2$.
More concretely, we show the following.
There exist twistor spaces on $4\mathbb{CP}^2$ such that 
(i) the anticanonical system is 4-dimensional as a linear system, and 
the image of the associated rational map is a scroll $Y$ of 2-planes in $\mathbb{CP}^4$ over a conic,
(ii) there is an explicit and simple elimination of the indeterminacy locus of the anticanonical map, whose  resulting morphisms is a
generically 2 to 1 covering map onto the scroll $Y$,
(iii) the branch divisor of the last covering, which will be denoted by $B$ throughout this paper, is  an intersection of \,$Y$\! with a  quartic hypersurface in $\mathbb{CP}^4$,
(iv) if we take homogeneous coordinates on $\mathbb{CP}^4$ such that the scroll $Y$ is defined by 
 $z_0^2=z_1 z_2$, 
then the quartic hypersurface is defined by the equation
\begin{align}\label{eqn:B}
z_0z_3z_4 f(z_0,z_1,z_2,z_3,z_4)=Q(z_0,z_1,z_2,z_3,z_4)^2,
\end{align}
 where $f$ and $Q$ are linear and quadratic polynomials with real coefficients respectively.
The double covering structure and  similarity of  the equations \eqref{eqn:B} with \eqref{KKP}  would justify to call these twistor spaces  a direct generalization of those by Kreussler-Kurke-Poon.

The main tool of the present investigation is the anticanonical 
system of the twistor spaces.
In Section \ref{s:1} we  start by constructing 
a rational surface $S$ which will be contained in the twistor space $Z$ on $4\mathbb{CP}^2$
as a real half-anticanonical divisor, and then 
clarify the structure of bi-anticanonical system on $S$.
Next in Section \ref{s:anticanonical} we study the structure of the anticanonical
map of the twistor spaces in detail.
In Section \ref{ss:acs} we show that the anticanonical map induces a rational map to $\mathbb{CP}^4$, whose image is a scroll $Y$ of planes over a conic, and give an explicit elimination of the 
indeterminacy locus, obtaining a degree 2 morphism $Z_1\to Y$.
Next in Section \ref{ss:modify}
we analyze the structure of the anticanonical map more
in detail, and 
by applying some explicit blowups and blowdowns we modify the degree 2 morphism $Z_1\to Y$ to another morphism $Z_4\to Y$ so as to have no divisor to be contracted. 
Up to contraction of curves, this gives the Stein factorization
of the morphism $Z_1\to Y$. 
Although the modifications therein are a little bit complicated,
they are rather natural in light of the structure of the anticanonical system, and indispensable for obtaining explicit construction 
of the twistor spaces. 
We also obtain a key technical result that the branch divisor of $Z_1\to Y$ is a cut of 
$Y$ by a quartic hypersurface.

In Section \ref{s:defeq}, we explicitly determine a defining  equation of the branch divisor of the double covering.
For this, in Section \ref{ss:dc} we  find 5 hyperplanes in $\mathbb{CP}^4$
such that the intersection of the branch divisor with the hyperplanes becomes double curves, i.e.\,a curve of multiplicity 2.
These double curves are analogous to the above 4 conics appeared in 
the case of $3\mathbb{CP}^2$, but
to understand how they  intersect each other requires some effort.
In addition, in the present case, finding all the double curves is not so easy, since not all of the double curves are obtained as an image of twistor lines as in the case of $3\mathbb{CP}^2$.
In Section \ref{ss:quadric} we show that the 5 double curves are
contained in a quadratic hypersurface in $\mathbb{CP}^4$, and that 
such a hyperquadric is unique up to the defining equation of 
the scroll $Y$.
In Section \ref{ss:defeq} we prove the main result which determines 
the defining equation of the quartic hypersurface (Theorem \ref{thm:B}). 
The equation includes not only the quadratic polynomial obtained in Section \ref{ss:quadric} but also a linear polynomial, 
which might look strange at first sight.
We  give an account for a geometric meaning of it.
 
In Section \ref{ss:sing} we investigate  singularities of 
the branch divisor $B$ of the double covering.
As above $B$ has 5 double curves, and at most of the intersection points of them, $B$ has (non-real) ordinary double points.
This is totally parallel to the case of $3\mathbb{CP}^2$ explained at the beginning.
But in the present case there are exactly 2 special intersection points, at which $B$ has 
$A_3$-singularities.
Besides these ordinary double points and $A_3$-singular points,
we show that $B$ has other isolated singularities and determine their basic invariants (Theorem \ref{thm:6sing}).
The result  means that in general $B$ has extra 6 ordinary double points.
For obtaining this result,
as in Kreussler \cite{Kr89} and Kreussler-Kurke \cite{KK92} in the case of $3\mathbb{CP}^2$,
we compute the Euler numbers of the relevant spaces, especially the branch divisor $B$.
We also note that the concrete modification of the anticanonical map 
obtained in Section \ref{ss:modify} is crucial for determining 
the invariants of the singularities.

In Section \ref{ss:moduli} we compute dimension of the moduli
space of the present twistor spaces. 
We first compute the dimension by counting the number of effective parameters (coefficients) involved in the defining
quartic polynomials.
Next we show that some cohomology group of the twistor space
an be regarded as a tangent space of  
the moduli space, and see that it coincides with the 
dimension obtained by counting the number of effective parameters.
This implies a completeness of our description obtained in Section \ref{s:defeq}.
In Section \ref{ss:gen} 
 we discuss genericity of our twistor spaces 
among all Moishezon twistor spaces on $4\mathbb{CP}^2$,
indicated in the title of this paper.
The genericity implies a kind of density in the moduli space of Moishezon twistor spaces on $4\mathbb{CP}^2$.
In the course we also give a rough classification of Moishezon twistor spaces on $4\mathbb{CP}^2$ under some genericity assumption. 

In Appendix we show that an inverse of a non-standard contraction map employed 
in Section \ref{ss:modify} can be realized by an embedded blowup with a non-singular center
in the ambient space, and  
point out in a concrete form that the modification can be
regarded as a singular version of the well-known operation of Hironaka
\cite{Hi60} for constructing non-projective Moishezon threefolds. 

We should also  mention a relationship between this work and our previous paper \cite{Hon-II}.
As we explained in the beginning, 
in the case of $3\mathbb{CP}^2$, the branch quartic divisor of the double covering becomes of the form \eqref{KKP} under a genericity assumption.
In non-generic cases, as showed by Kreussler-Kurke \cite{KK92}, the branch divisor becomes similar but more degenerate form, and in the most degenerate situation, the branch divisor has a $\mathbb C^*$-action.
As a consequence, in that case the twistor spaces have a $\mathbb C^*$-action.
The twistor spaces studied in \cite{Hon-II} is a generalization of these twistor spaces (with $\mathbb C^*$-action) to the case of $n\mathbb{CP}^2$, $n>3$.
In this respect we remark that the twistor spaces in this paper is obtained as a deformation of the twistor spaces studied in \cite{Hon-II} in the case of $4\mathbb{CP}^2$. 
It is very natural to expect that we can obtain a further generalization 
to the case of $n\mathbb{CP}^2$, $n>4$.
We hope to discuss this attractive topic in a future paper.

\vsp
\noindent
{\bf Acknowledgement.} I would like to thank Masaharu Ishikawa for kindly answering a question about deformation of singularities, which helps me to be confident about  computations of the Euler number of the singular branch surface. 

\noindent
{\bf Notations and Conventions.}
The natural square root of the anticanonical bundle of 
a twistor space is denoted by $F$.
If two varieties $X_1$ and $X_2$ are birational under some blowup or blowdown, and if $Y$ is a subvariety of $X_1$, then we often use the same symbol $Y$ to mean the strict transform or birational image of $Y$ under the blowup or blowdown, as far as it makes sense. 
For a line bundle $L$ and a non-zero section $s$
of $L$, $(s)$ means the zero divisor of $s$.
For a linear subspace $V\subset H^0(L)$, we denote by $|V|$ to mean the linear system $\{(s)\set s\in H^0(L),\,s\neq0\}$.
$\Bs \,|V|$ means the base locus of $|V|$.
We mean $\dim |V| = \dim V -1$ and $h^i(L)=\dim H^i(L)$.

\section{A construction of rational surfaces and their bi-anticanonical system}\label{s:1}

\begin{figure}
\includegraphics{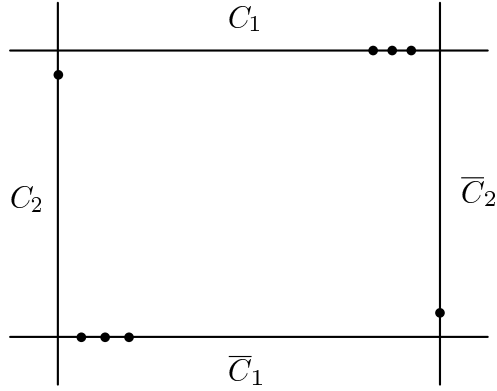}
\caption{The 8 points to be blown up, giving the surface $S$.}
\label{fig:square}
\end{figure}

We are going to investigate  twistor spaces on $4\mathbb{CP}^2$ which contain 
a particular type of non-singular rational surface $S$ as a real member of $|F|$.
In this section we first construct the surface $S$ as a blowup of $\mathbb{CP}^1\times\mathbb{CP}^1$, and then study the bi-anticanonical system on it.
For this, we define the line bundles 
$\mathscr O(1,0)$ and $\mathscr O(0,1)$ 
on $\mathbb{CP}^1\times\mathbb{CP}^1$ as the pullback of $\mathscr O_{\mathbb{CP}^1}(1)$ by 
the projection to the first and second factors respectively.
We simply call members of the linear system $|\mathscr O
(m,n)|$ as $(m,n)$-curves.
We define a real structure on  $\mathbb{CP}^1\times\mathbb{CP}^1$ as a product of the complex conjugation and the antipodal map.
Next take any non-real $(1,0)$-curves $C_1$, any $(0,1)$-curve $C_2$, any distinct  3 points on $C_1\backslash(C_2\cup\ol C_2)$, and 1 point on $C_2\backslash(C_1\cup\ol C_1)$.
By taking the images under the real structure, we obtain distinct 8 points in total (see Figure \ref{fig:square}).
Let $\epsilon:S\to \mathbb{CP}^1\times\mathbb{CP}^1$ be the blowup at these 8 points.
This surface $S$ has a natural real structure induced from that on $\mathbb{CP}^1\times\mathbb{CP}^1$, and also has an anticanonical curve 
$$
C:=C_1+C_2+\ol C_1+\ol C_2,
$$
where this time $C_i$ means the strict transform of the original $C_i$,
so that $C_1^2=\ol C_1^2=-3$ and $C_2^2 = \ol C_2^2= -1$.
The identity component of the holomorphic automorphism group of $S$ 
is trivial.
In the sequel we investigate the anticanonical and bi-anticanonical systems on $S$.

\begin{proposition}\label{prop:S1}
(i) $\dim |K_S^{-1}|=0$, so that $C$ is the unique anticanonical curve on $S$,
(ii) $\dim |2K_S^{-1}|=2$,  $\Bs\,|2K_S^{-1}|=C_1\cup\ol C_1$, and $\,\Bs\,|2K_S\inv-C_1 - \ol C_1|=\emptyset$. 
\end{proposition}

\proof
(i) is immediate.
For (ii),  $C C_1=C\ol C_1=-1$ means the inclusion $C_1\cup\ol C_1\subset\Bs\, |2K_S\inv|$.
Riemann-Roch formula and the rationality of $S$ imply $\chi(K_S \inv)=1$.
These mean $H^1(K_S\inv)=0$.
Therefore  
restricting  $2K_S\inv -  C_1 - \ol C_1 $ to $C_2\cup \ol C_2$, we  obtain the exact sequence
\begin{align}\label{es1}
0 \lra H^0 (K_S\inv) \lra H^0 (2K_S\inv -  C_1 - \ol C_1 )
\lra H^0(\mathscr O_{C_2})\oplus H^0(\mathscr O_{\ol C_2}) \lra 0.
\end{align}
This implies  $(C_2\cup \ol C_2)\cap \Bs \, |2K_S\inv - C_1 - \ol C_1 |=\emptyset$, and also $h^0 (2K_S\inv -  C_1 - \ol C_1 )=3 $.
On the other hand by restricting the same system to $C_1\cup\ol C_1$, we obtain the exact sequence 
\begin{align}\label{es2}
0 \lra H^0(\mathscr O_{S}(2C_2+2\ol C_2)) \lra H^0(2K_S\inv-C_1-\ol C_1) 
\stackrel{r}{\lra} H^0(\mathscr O_{C_1}(1))\oplus H^0(\mathscr O_{\ol C_1}(1)).
\end{align}
As $h^0(\mathscr O_{S}(2C_2+2\ol C_2))=1$ clearly, the image of the restriction map $r$ in \eqref{es2} is 2-dimensional.
Projecting this to $H^0(\mathscr O_{C_1}(1))$ gives either a 1-dimensional subspace or $H^0(\mathscr O_{C_1}(1))$ itself.
In order to prove $\Bs\,|2K_S\inv - C_1 - \ol C_1| = \emptyset$, it is enough to exclude the former possibility.
Suppose it is the case.
Let $p\in C_1$ be the zero point of a generator of the 1-dimensional subspace. 
Then we have $\Bs\,|2K_S\inv - C_1- \ol C_1| = \{ p, \ol p \}$ and $p\not\in C_2\cup\ol C_2$ as we have already seen.
Let $S'\to S$ be the blowup at $p$ and $\ol p$, and $C'_1$ and $\ol C'_1$ the strict transforms of $C_1$ and $\ol C_1$ respectively.
As $(2K_S\inv - C_1 - \ol C_1)^2 = 2$, for any member $D\in |2K_S\inv - C_1 - \ol C_1|$ with 
$D\neq C_1+2C_2+ \ol C_1 + 2\ol C_2$,
we have $C\cap D=\{p,\ol p\}$, and the intersections are transversal.
This means that the system  $| C'_1+ 2C_2+\ol C'_1+ 2\ol C_2 |$ on $S'$ is base point free.
Then as $C'_1( C'_1+ 2C_2+\ol C'_1+ 2\ol C_2 ) = \ol C'_1( C'_1+ 2C_2+\ol C'_1+ 2\ol C_2 ) = 0$,
the morphism associated to the system contracts $C'_1$ and $\ol C'_1$ to points.
On the other hand, by \eqref{es1}, $C_2$ and $\ol C_2$ are mapped to mutually different points by the same morphism.
This is a contradiction because $C_1\cap C_2\neq\emptyset$, $C_1\cap \ol C_2\neq\emptyset$,
and $C_1$ is connected.
\proofend

\vsp
Let $\phi:S\to \CP^2$ be the morphism associated to the system $|2K_S \inv | \simeq |2K_S\inv - C_1 - \ol C_1|$.
Since $\Bs\,|2K_S^{-1}-C_1-\ol C_1| = \emptyset$ and 
$(2K_S^{-1}-C_1-\ol C_1)C_2 = (2K_S^{-1}-C_1-\ol C_1)\ol C_2 =0$,
$\phi$ factors as $S\to \ol S\to \mathbb{CP}^2$, where 
$S\to \ol S$ denotes the blowdown of $C_2$ and $\ol C_2$.

\begin{proposition}\label{prop:2to1}
The morphism $\phi$ is generically 2 to 1, and the branch divisor is a quartic curve.
Further, the images $\phi(C_1)$ and $\phi(\ol C_1)$ are the same line,
and $\phi(C_2)$ and $\phi ( \ol C_2 )$ are different 2 points on the line.
\end{proposition}

\proof
The morphism $\phi$ is surjective since $(2K_S\inv - C_1 - \ol C_1)^2 = 2>0$.
This also means it is generically 2 to 1.
Further, a general member $D$ of $|2K_S\inv - C_1 - \ol C_1 |$, which is an irreducible non-singular curve by Bertini's theorem,
is an elliptic curve because $D(D+K_S)=0$.
Hence the branch curve of $\phi$ must be of degree 4.
For the images $\phi(C_1)$ and $\phi(\ol C_1)$, as 
in the proof of Proposition \ref{prop:S1},
we have the exact sequence \eqref{es2}, and the image of $r$
projects isomorphically to $H^0(\mathscr O_{C_1}(1))$ and $H^0(\mathscr O_{\ol C_1}(1))$. 
Hence $\phi(C_1)$ and $\phi(\ol C_1)$ are lines.
Further these lines are identical, since
they are precisely the 2-dimensional images into the dual space
$\mathbb P H^0(2K_S^{-1}-C_1-\ol C_1)^*$ by the 
dual map of $r$.
Also the last claim for $\phi(C_2)$ and $\phi(\ol C_2)$ is clear from the exact sequence \eqref{es1}, $C_1\cap C_2\neq\emptyset$ and $C_1\cup \ol C_2\neq\emptyset$.
\proofend

\vsp
The following property of the branch quartic will also be needed later.
\begin{proposition}\label{prop:S3}
If we denote the line and  the pair of points on it by $l:=\phi(C_1)=\phi (\ol C_1)$ and
$p_2:=\phi(C_2)$ and $\ol p_2=\phi(\ol C_2)$ respectively,
then $p_2$ and $\ol p_2$ are smooth points of the branch quartic curve,
and the curve is tangent to $l$ at these 2 points.
(Hence $l$ is a bitangent of the branch curve.)
\end{proposition}

\proof
Let $\beta$ be the branch curve of $\phi$.
As $C_1+2C_2+ \ol C_1 + 2\ol C_2\in |2K_S\inv - C_1 - \ol C_1|$, 
there is a line $l'$ such that  $\phi^{-1}(l')=C_1+2C_2+ \ol C_1 + 2\ol C_2$.
But as $\phi(C_1)=l$, we obtain $l'=l$.
Since $\phi(C_2)=p_2$ and both $C_2\cap C_1$ and $C_2\cap \ol C_1$ are non-empty,
it follows that $p_2\in \beta$, so that $\ol p_2\in\beta$ by the real structure.
Further, after the blowdown $S\to \ol S$, the curve $C_1\cup\ol C_1$ is of course locally reducible at the images  of 
$C_2$ and $\ol C_2$.
Hence we have  $\beta|_{l}=2p_2+2\ol p_2$ as divisors. 
This means that either $\beta$ has double points at $p_2$ and $\ol p_2$,
or otherwise $\beta$ is smooth at $p_2$ and $\ol p_2$ and is tangent to $l$ at these 2 points.
But in the former case  by smoothness of $S$ there have to be extra exceptional curves over $p_2$ and $\ol p_2$, which contradicts $\phi^{-1}(l)= C_1+2C_2+ \ol C_1 + 2\ol C_2$.
Hence the claim follows.
\proofend

\section{Analysis of the anticanonical map on the twistor spaces}
\label{s:anticanonical}
\subsection{The anticanonical map of the twistor spaces}
\label{ss:acs}
Let $S$ be the rational surface equipped with the real structure  constructed in the previous section, and $C=C_1+C_2+\ol C_1+\ol C_2$ the unique 
anticanonical curve on $S$.
Let $Z$ be a twistor space on $4\CP^2$ and suppose that $Z$ contains $S$ as a 
real member of $|F|$.
The following property of $|F|$ is immediate to see and we omit a proof.

\begin{proposition}\label{prop:F}
The system $|F|$ satisfies the following:
(i) $\dim |F|=1$, (ii) $\Bs\,|F| = C$,
(iii) the number of reducible members of $|F|$ is two, and
both of the members are real.
\end{proposition}

\begin{figure}
\includegraphics{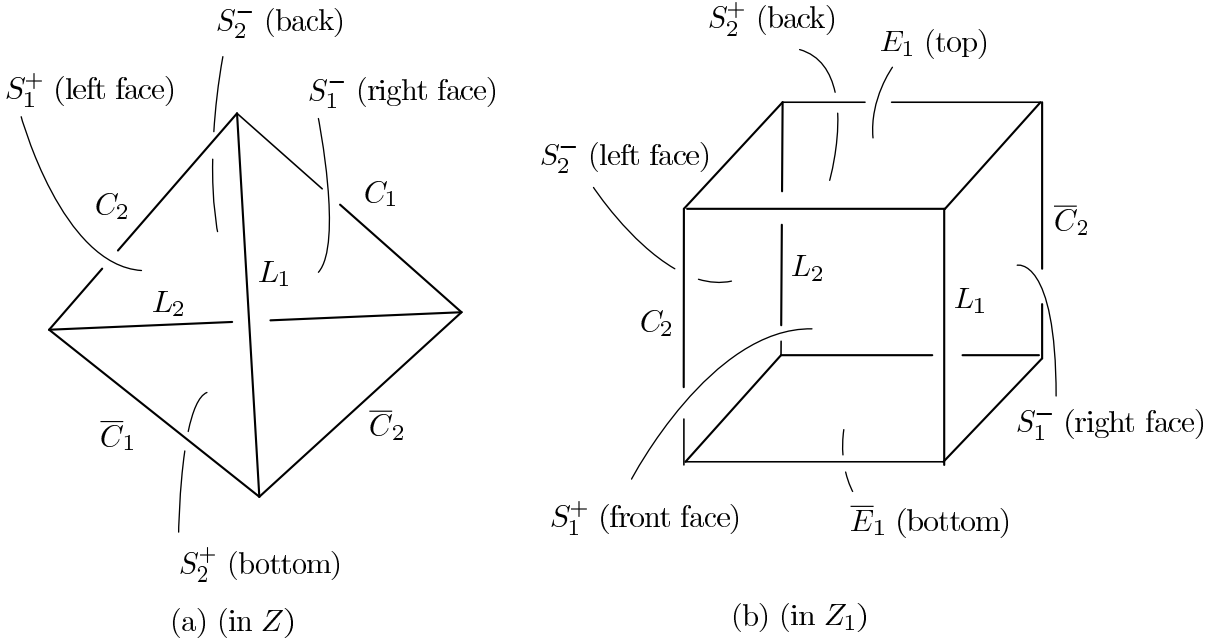}
\caption{a tetrahedron in $Z$ and a cube in $Z_1$}
\label{fig:trans0}
\end{figure}
We note that it readily follows from (i) and (ii) that a general member $S'$ of the pencil $|F|$ is also
obtained from $\mathbb{CP}^1\times\mathbb{CP}^1$ by blowing up
8 points arranged as in Figure \ref{fig:square},
where the positions of the 8 points are not identical
to the original ones.
We denote the 2 reducible members of $|F|$ by 
\begin{align}\label{reduciblemembers}
S_1=S_1^++S_1^-\quad{\rm{and}}\quad S_2=S_2^++S_2^-,
\end{align}
 where we make distinction between $S_i^+$ and $S_i^-$ by declaring that 
$S_1^+$ and $S_2^+$ contains the component $\ol C_1$.
We denote $L_1:=S_1^+\cap S_1^-$ and $L_2:=S_2^+\cap S_2^-$, both of which are twistor lines by \cite[\S 1]{P92}.
Then these divisors and curves  form
a tetrahedron as illustrated
in Figure \ref{fig:trans0}, (a).
These will be significant for our analysis of the anticanonical system on the twistor spaces.
We  show the the following basic properties of the anticanonical system. 
Note that (iv) means that $Z$ is Moishezon.

\begin{proposition}\label{prop:2F}
The anticanonical system $|2F|=|K_Z\inv|$ of the twistor space $Z$ satisfies the following:
(i) $\dim |2F|=4$, (ii) $\Bs\, |2F| = C_1 \cup \ol C_1$,
(iii) if $\mu_1: Z_1\to Z$ denotes the blowup at $C_1\cup \ol C_1$,  $E_1\cup\ol E_1$ the exceptional divisor, and $\mathscr L_1:= \mu_1^*(2F)-E_1-\ol E_1$, then
$\Bs \, |\mathscr L_1| = \emptyset$,
(iv) if $\Phi_1$ denotes the morphism associated to $|\mathscr L_1|$,
then the image ${\Phi}_1(Z_1)$ is a  scroll of 2-planes over a conic, and
 the morphism ${\Phi}_1$ is generically 2 to 1 over the scroll.
\end{proposition}

By the blowup $\mu_1:Z_1\to Z$, the tetrahedron in $Z$ is 
transformed to be a cubic in $Z_1$ as in Figure \ref{fig:trans0}, (b).

\noindent
{\em Proof of Proposition \ref{prop:2F}.}
(i) is immediate from Proposition \ref{prop:S1} (ii), Proposition \ref{prop:F} (i), and the exact sequence
\begin{align}\label{ses:basic1}
0 \lra H^0 (F) \lra H^0 (2F) \lra H^0(2K_S\inv) \lra 0,
\end{align}
where the last zero is a consequence of $h^0(F)=2$ and the Riemann-Roch formula applied to $F$.
The claim (ii)  also follows from this exact sequence and Proposition \ref{prop:S1} (ii). 

For (iii), let $\tilde S$ be the strict transform of $S$.
Then $\tilde S$ is biholomorphic to $S$, and $\tilde S\in |\mu_1^*F-E_1-\ol E_1|$.
Hence we have an exact sequence
\begin{align}\label{es3}
0 \lra \mathscr \mu_1^*F \lra 
\mu_1^*(2F)-E_1-\ol E_1 \lra \mu_1^*(2F)-E_1-\ol E_1|_{\tilde S} \lra 0. 
\end{align}
Since $H^1(\mu_1^*F)\simeq H^1(F)=0$, we obtain that
the restriction map 
$H^0(\mu_1^*(2F)-E_1-\ol E_1) \to H^0(\mu_1^*(2F)-E_1-\ol E_1|_{\tilde S})$ is surjective.
Further,
as $(\mu_1^*F)|_{\tilde S}\simeq F|_S\simeq K_S \inv$ under the biholomorphism $\tilde S\simeq S$,
we have $\mu_1^*(2F)|_{\tilde S}\simeq 2K_S \inv$.
Further, $E_1|_{\tilde S}\simeq \mathscr O_S(C_1)$.
Hence we obtain an isomorphism $\mu_1^*(2F)-E_1-\ol E_1|_{\tilde S} \simeq 2K_S\inv - C_1 - \ol C_1$.
Therefore by the third claim of Proposition \ref{prop:S1} (ii), we obtain $\Bs \, |\mu_1 ^*(2F)-E_1-\ol E_1| = \emptyset$.

Let $\Phi:Z\to \CP^4$ be the rational map associated to the anticanonical system $| 2F |$, so that ${\Phi_1}=\Phi\circ\mu_1$.
For (iv)  it is enough to show that $\Phi(Z)$ is the 3-dimensional scroll as in the statement, 
and the rational map $\Phi:Z\to \Phi(Z)$ is generically 2 to 1.
Let $S^2 H^0(F)$ be the subspace of $H^0(2F)$ generated by all sections of the form $s_1s_2$ where $s_i\in H^0(F)$. 
This is a 3-dimensional subspace.
Then we have the following left commutative diagram of rational maps:

\begin{equation}\label{cd1}
\xymatrix{
   Z \ar@{->}[r]^{{\Phi}}  \ar@{->}[rd]_{{f}}  & \mathbb{CP}^{4} \ar@{->}[d]^{\pi} \\
  & \mathbb{CP}^{2}\\
}
\qquad\qquad
\xymatrix{
  Z \ar@{->}[r]^{{\Phi}}  \ar@{->}[rd]_{{f}}  & Y \ar@{->}[d]^{\pi} \\
  & \Lambda  
}
 \end{equation}
where $\pi$ is the linear projection induced by the inclusion $S^2H^0(F)\subset H^0(2F)$
and $f$ is the rational map associated to the subsystem $|S^2H^0(F)|$.
Clearly the image $f(Z)$ is a conic, for which we denote by $\Lambda$.
Hence writing $Y:=\pi^{-1}(\Lambda)$, $Y$ is exactly the  scroll as in the statement of (iv),
and we  obtain the right commutative diagram in \eqref{cd1}.
We have to show that $\Phi:Z\to Y$ is surjective and generically 2 to 1.
For these, we note that by the definition of $f$, for  any $\lambda\in \Lambda$, $f^{-1}(\lambda)$ belongs to the pencil $|F|$.
Then by \eqref{ses:basic1} for {\em any} non-singular member $S\in |F|$,  the restriction $\Phi|_S$ is exactly 
the rational map associated to the system $|2 K_S \inv |$,
where the target space is  the fiber plane  $\pi^{-1}(\lambda)$.
By Proposition \ref{prop:2to1}, this means that $\Phi|_{f^{-1}(\lambda)}:f^{-1}(\lambda)\to \pi^{-1}(\lambda)$ is surjective as far as $f^{-1}(\lambda)$ is non-singular.
Therefore $\Phi$ itself is also surjective to $Y$.
Now the final claim (2 to 1 over $Y$) is immediate from these considerations and Proposition \ref{prop:2to1}.
\proofend

\vsp
Thus the anticanonical map is 2 to 1 over the scroll $Y$.
Further, by the above argument and Proposition \ref{prop:2to1}, 
the branch locus of the 2 to 1 map has degree 4 on  the planes $f^{-1}(\lambda)$,
from which one might find similarity with those in the case of $3\mathbb{CP}^2$ \cite{KK92, P92}.
But we are yet far from the goal.
In the next subsection we shall investigate  structure of the anticanonical map more closely.

\subsection{Modification of the anticanonical map}\label{ss:modify}
We use the notations $\Phi$, $\Lambda$, $Y$, $f$ and $\pi$ given in the proof of the  proposition.
Further define $l$ to be the singular locus of the scroll $Y$.
$l$ is a line, and is exactly the  indeterminacy locus of the projection $\pi$.
This line plays an important role throughout this paper.
Note that for a hyperplane $H\subset\CP^4$, the intersection $Y\cap H$ splits to planes iff $H$ projects to a line in $\mathbb{CP}^2$ (in which the conic $\Lambda$ is contained), and otherwise $Y\cap H$ is a cone over $\Lambda$ whose vertex is the point $l\cap H$.
Further, in the former situation, $Y|_H$ is a double plane (i.e.\,a non-reduced plane of multiplicity 2) iff the line is  tangent to $\Lambda$. 
Let $\nu:\tilde Y\to Y$ be the blowup at $l$, and $\Sigma$ the exceptional divisor.
$\tilde Y$ is biholomorphic to the total space of the $\mathbb{CP}^2$-bundle 
$\mathbb P(\mathscr O(2)^{\oplus 2}\oplus \mathscr O)\to \Lambda$,
and $\Sigma$ is identified with the subbundle $\mathbb P(\mathscr O(2)^{\oplus 2})$, so that it is
biholomorphic to $\mathbb{CP}^1\times\mathbb{CP}^1$.
More invariantly, we have a natural isomorphism
$\Sigma\simeq l\times\Lambda$.
The composition $\tilde Y\to Y\to \Lambda$ is a morphism which is identified with the bundle projection, for which we denote by $\tilde{\pi}$.

In order to treat the branch locus of the degree 2 morphism $\Phi_1$ (see Proposition \ref{prop:2F}) properly, we consider a lifting problem of the morphism  $\Phi_1$ to $\tilde Y$.
By definition of $f$, the composition $Z_1\stackrel{\mu_1}{\to} Z\stackrel{f}{\to} \Lambda\subset\mathbb{CP}^2$ is the rational map associated to the 2-dimensional linear system
$\mu_1^*|S^2H^0(F)|$ on $Z_1$ (whose members are the total transforms of members of the linear system $|S^2H^0(F)|$).
This linear system on $Z_1$ is a subsystem of $|\mu_1^*(2F)|$.
But since the pencil $|F|$ has $C_1$ and $\ol C_1$ as  components of the base locus,
all the above total transforms contain the divisor $2E_1+2\ol E_1$.
Hence subtracting $E_1+\ol E_1$ and recalling $\mathscr L_1=\mu_1^*2F-E_1-\ol E_1$, 
we can regard linear system $\mu_1^*|S^2H^0(F)|$ as a subsystem of $|\mathscr L_1|$.
This 2-dimensional subsystem of $|\mathscr L_1|$  still has $E_1+\ol E_1$ as the fixed components, so by subtracting it we obtain a 2-dimensional  subsystem of $|\mathscr L_1-E_1-\ol E_1|=|\mu_1^*2F-2(E_1+\ol E_1)|$,
which is readily seen to be coincide with $|\mu_1^*2F-2(E_1+\ol E_1)|$ itself.
Hence the composition $Z_1\to Z\to \Lambda$ can be regarded as the rational map associated to $|\mu_1^*2F-2(E_1+\ol E_1)|$.
However, the curves $C_2$ and $\ol C_2$ are contained in Bs\,$|F|$
(Proposition \ref{prop:F} (ii)), and the strict transforms of these curves to $Z_1$ are exactly the base locus of $|\mu_1^*2F-2(E_1+\ol E_1)|$.
Thus the composition $Z_1\to Z\to \Lambda$ has the strict transforms of $C_2$ and 
$\ol C_2$  as its indeterminacy locus.
Therefore, the morphism $\Phi_1:Z_1\to Y$ cannot be lifted to  $Z_1\to \tilde Y$ as a morphism, because the composition with $\tilde{\pi}:\tilde Y\to\Lambda$ is not a morphism.

So let $\mu_2: Z_2\to Z_1$ be the blowup at $C_2\cup \ol C_2$,
and $E_2$ and $\ol E_2$ the exceptional divisors.
Here we are regarding $C_2$ and $\ol C_2$ as curves in $Z_1$.
Let $\Phi_2:=\Phi_1\circ\mu_2$, and $\mathscr L_2:=\mu_2^*\mathscr L_1$.
(This time we do not subtract $E_2+\ol E_2$ because $C_2$ and $\ol C_2$ are not base curves of  $|\mathscr L_1|$.)
Obviously $\Phi_2$ is the rational map associated to $|\mathscr L_2|$, and it is clearly a morphism.
Then we have the following:

\begin{proposition}
The morphism $\Phi_2:Z_2\to Y$ can be lifted to a morphism $\tilde{\Phi}_2: Z_2\to\tilde{Y}$.
Namely there is a morphism $\tilde{\Phi}_2:Z_2\to\tilde Y$ such that $\Phi_2$ factors as $Z_2
\stackrel{\tilde{\Phi}_2}{\to} \tilde Y \stackrel{\nu}{\to} Y$.
\end{proposition}

\proof
As in the above explanation, the composition $Z_2\to Z_1\to \Lambda
\subset\mathbb{CP}^2$ is the rational map associated to the system $|\mu_2^*\{\mu_1^*2F-2(E_1+\ol E_1)\}|$.
In the same way for the identification between a conic and a line on a plane by means of
the projection from a point, once we fix any 
non-zero element of $ H^0(\mu_1^*F-E_1-\ol E_1)$,
by taking a product with it,
members of the system
 $|\mu_1^*2F-2(E_1+\ol E_1)|$ can be identified with those of
the pencil $|\mu_1^*F-(E_1+\ol E_1)|$,
and the rational map associated to $|\mu_1^*2F-2(E_1+\ol E_1)|$
is identified with the rational map associated to
$|\mu_1^*F-(E_1+\ol E_1)|$.
Therefore the composition $Z_2\to Z_1\to \Lambda
\subset\mathbb{CP}^2$ can be identified with the rational map associated to the pencil $|\mu_2^*(\mu_1^*F-E_1-\ol E_1)|$.
This pencil has $E_2+ \ol E_2$ as the fixed component,
and if we subtract this, the pencil becomes free, since 
$|F|_S|=|K_S^{-1}|$ and $|K_S^{-1}|$ consists of a single member $C_1+\ol C_1+C_2+\ol C_2$ (a reduced curve) by Proposition \ref{prop:S1} (i).
Hence the composition $Z_2\to Z_1\to \Lambda$ has no point of indeterminacy.
We write $f_2$ for this morphism.
Thus we are in the following left situation :
\begin{equation}\label{cd2}
\xymatrix{
   Z_2 \ar@{->}[d]_{f_2}  \ar@{->}[rd]  & 
   \tilde Y \ar@{->}[d]^{\nu} \ar@{->}[ld] \\
 \Lambda & Y \ar@{->}[l]^{\pi}
}
\qquad\qquad
\xymatrix{
   Z_2 \ar@{->}[d]_{f_2}  \ar@{->}[r]^{\tilde{\Phi}_2}  \ar@{->}[rd]^{\Phi_2}   & 
   \tilde Y \ar@{->}[d]^{\nu}  \\
 \Lambda & Y \ar@{->}[l]^{\pi}
}
 \end{equation}
 where all maps except $\pi$ are morphisms, and  the 2 triangles are commutative. (The map $Z_2\to Y$ is $\Phi_2$ and the map $\tilde Y\to \Lambda$ is $\tilde{\pi}$.)
For lifting $\Phi_2$ to $\tilde Y$, we need to assign a point of $\tilde Y$ for each point of $Z_2$.
For this, recall that $\nu$ is isomorphic outside $\nu^{-1}(l)=\Sigma$.
For $z\in Z_2\backslash\Phi^{-1}_2(l)$, of course, we assign 
the point $\nu^{-1}(\Phi_2(z))$.
For $z\in \Phi_2^{-1}(l)$ define $\lambda:=f_2(z)$ and $y:=\Phi_2(z)\in l$.
The inverse image $\nu^{-1}(y)$
is a fiber of the projection $\Sigma\to l$.
In accordance with the natural isomorphism $\Sigma\simeq l\times \Lambda$,
this fiber and the fiber $\tilde{\pi}^{-1}(\lambda)$ intersect at a unique point.
Let $\tilde y\in\tilde Y$ be this point, and we assign $\tilde y$ to $z$.
Define  $\tilde{\Phi}_2:Z_2\to\tilde Y$  to be the map thus obtained.
Then $\tilde{\Phi}_2$ is clearly continuous and is a lift of $\Phi_2$.
As $\tilde{\Phi}_2$ is holomorphic on the complement of the analytic subset $\Phi_2^{-1}(l)$, Riemann's extension theorem means that  it is automatically holomorphic on the whole of $Z_2$.
Thus we get the  situation right in \eqref{cd2} and obtained the desired lift $\tilde{\Phi}_2$.
\proofend

\vspace{3mm}
Since the lift $\tilde{\Phi}_2:Z_2\to \tilde Y$ is a degree 2 morphism between non-singular spaces, we can speak about its branch divisor.
Namely we first define the ramification divisor $R$ on $Z_2$  as a zero divisor of a natural section 
(defined by the Jacobian) of the line bundle
$K_{Z_2} - \tilde{\Phi}_2^* K_{\tilde Y}$, and then let the branch divisor $\tilde B$ to be
the image $\tilde{\Phi}_2(R)$, which is necessarily a divisor.
We can determine the cohomology class of this divisor as follows:

\begin{proposition}\label{branchdivisor1}
Let $\tilde B$ be the branch divisor of the lift $\, \tilde{\Phi}_2: Z_2\to \tilde Y$
as above.
Then $\tilde B\in |\mathscr O_{\tilde Y}(4)|$, where $\mathscr O_{\tilde Y}(1):=
\nu^*\mathscr O_Y(1) = \nu^*\mathscr O_{\mathbb{CP}^4}(1)|_Y$.
\end{proposition}

\proof
As before let $\Sigma$ be the exceptional divisor of the blowup $\tilde Y\to Y$,
and let $\mathfrak f$ be the cohomology class of the fiber class $\tilde{\pi}^*\mathscr O_{\Lambda}(1)$.
The cohomology group $H^2(\tilde Y,\mathbb Z)$ is a free $\mathbb Z$-module generated by $\Sigma$ and $\mathfrak f$.
$\Sigma$ is isomorphic to $l\times\Lambda$,
and the restriction  $\nu|_{\Sigma}$
can be identified with the projection to $l$.
Define $(0,1)$ to be the bidegree of a fiber of this projection.
Then  the normal bundle is $N_{\Sigma/\tilde Y}\simeq\mathscr O(-2,1)$,
while $\mathfrak f$ is restricted to the class $(1,0)$.
From these we can readily deduce that the restriction map $H^2(\tilde Y,\mathbb Z)\to H^2(\Sigma,\mathbb Z)$ is isomorphic.

Let $H\subset\mathbb{CP}^4$ be any hyperplane containing the line $l$.
Then since $\Lambda$ is a conic, we have $\nu^{-1}(H) = \Sigma + 2\mathfrak f$, which means
\begin{align}\label{hs22}
\nu^*\mathscr O_Y(1) = \Sigma + 2\mathfrak f \hsp {\rm{in}} \hsp H^2(\tilde Y,\mathbb Z).
\end{align}
Then by using the above explicit form of the restriction map, we obtain 
\begin{align}
\mathscr O_{\tilde Y}(1)|_{\Sigma} \simeq
\Sigma|_{\Sigma} + 2\mathfrak f|_{\Sigma} = 
\mathscr O(-2,1) + \mathscr O(2,0)= \mathscr O(0,1).
\end{align}
Hence in order to prove $\tilde B\in |\mathscr O_{\tilde Y}(4)|$, it suffices to
show $\tilde B|_{\Sigma}\in|\mathscr O(0,4)|$.

In order to obtain the restriction 
$\tilde B|_{\Sigma}$, for each $\lambda\in \Lambda$ we write $S_{\lambda}:=f_2^{-1}(\lambda)$, which is the strict transform of a member of the pencil $|F|$.
Then it is not difficult to see that the restriction $\tilde{\Phi}_2|_{S_{\lambda}}:S_{\lambda}\to \tilde{\pi}^{-1}({\lambda})=\mathbb{CP}^2$ is naturally identified with the restriction of the original restriction $\Phi|_{S_{\lambda}}:S_{\lambda}\to\mathbb{CP}^2$.
If $S_{\lambda}$ is non-singular (which is the case for almost all $\lambda$), the last restriction $\Phi|_{S_{\lambda}}$ is exactly the bi-anticanonical map of $S_{\lambda}$.
Hence by Propositions \ref{prop:2to1} and \ref{prop:S3}, the branch curve is a  quartic curve which is tangent to the line $l$ at 2 points.
The last 2 points  are independent of the choice of $\lambda$, 
since  they are exactly $p_2=\Phi(C_2)$ and $\ol p_2=\Phi(\ol C_2)$.
Therefore $\tilde B|_{\Sigma}$ contains the fibers of $\Sigma\to l$  
(the restriction of $\tilde Y\to Y$ to $\Sigma$) over the points $p_2$ and $\ol p_2$ by multiplicity 2 respectively.
Moreover since the intersection of the line $l$ with the branch quartic curve  of  $S_{\lambda}\to \mathbb{CP}^2$ consists of the 2 points $p_2$ and $\ol p_2$, if $\tilde B|_{\Sigma}$  contains an irreducible curve of bidegree $(k,l)$ with $k>0$, then we have $(k,l)=(1,0)$.
But this cannot occur since the original $\Phi$ does not contain $l$ as a branch locus.
Thus we have obtained  $\tilde B|_{\Sigma}\in |\mathscr O(0,4)|$.
\proofend

\begin{figure}
\includegraphics{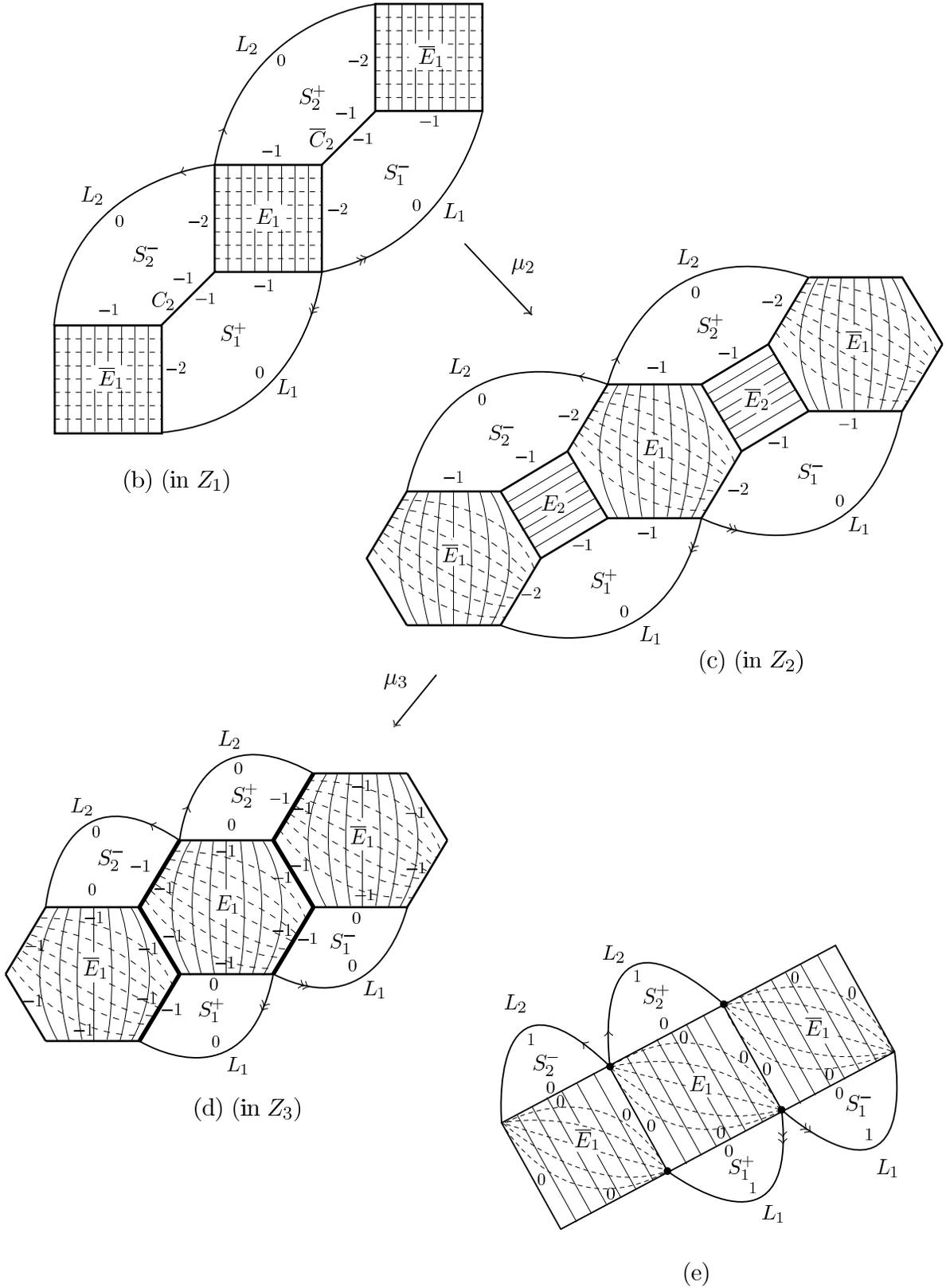}
\caption{The transformations from $Z_1$ to $Z_3$.
The picture (b) is identical to (b) in Figure \ref{fig:square}; the present (b) is obtained from 
the original (b) by just cutting out (just for presentation) along the 2 twistor lines
$L_1$ and $L_2$.
So
in each of (b), (c) and (d), the two $L_1$-s are identified in the direction indicated by 
the arrows, and the same for $L_2$-s.
(e) is obtained from (d) by contracting four $(-1,-1)$-curves, and will be used later.}
\label{fig:trans}
\end{figure}

\vspace{3mm}
Define a divisor $B$ on $Y$ by $B:=\nu(\tilde B)$.
By Proposition \ref{branchdivisor1}, $B\in |\mathscr O_Y(4)|$,
and since $B=\nu(\tilde{\Phi}_2(R))=\Phi_2(R)$, the morphism $\Phi_2:Z_2\to Y$ is a generically 2 to 1 covering with branch $B$.
(The former implies that $B$ is a cut of $Y$ by a quartic hypersurface. 
We will explicitly obtain a defining equation of this hypersurface in the next section.)
However $\Phi_2$ is {\em not}   a finite map but contracts the divisors $E_1,\ol E_1,E_2$ and $\ol E_2$ as we  see next until completing as Proposition \ref{prop:factor1}.
For this, we first notice that the exceptional divisors $E_2$ and $\ol E_2$ of $\mu_2:Z_2\to Z_1$ are isomorphic to $\mathbb{CP}^1\times\mathbb{CP}^1$ and the normal bundles satisfy
\begin{align}
N_{E_2/Z_2}\simeq \mathscr O(-1,-1),\quad
N_{\ol E_2/Z_2}\simeq \mathscr O(-1,-1).
\end{align} 
(See Figure \ref{fig:trans} (c).)
Therefore $E_2$ and $\ol E_2$ can also be blowdown along the projection different from the original $E_2\to C_2$ and $\ol E_2\to \ol C_2$.
Let $\mu_3:Z_2\to Z_3$ be this blowdown.
(See (c) $\to$ (d) in Figure \ref{fig:trans}.)
$Z_3$ is still non-singular.
The birational transformation from $Z_1$ to $Z_3$ is exactly Atiyah's flop at $C_2$ and $\ol C_2$.
The divisors $E_1$ and $\ol E_1$ in $Z_1$ are also isomorphic to $\mathbb{CP}^1\times\mathbb{CP}^1$, and they are respectively blown up at 2 points through $\mu_2$. 
We use the same letters $E_1$ and $\ol E_1$ to mean these divisors in $Z_2$.
These divisors are not affected by the blowdown $\mu_3$, and we still denote by $E_1$ and $\ol E_1$ for their images  in $Z_3$, as displayed in Figure \ref{fig:trans} (d).

Next we show that these 2 divisors $E_1$ and $\ol E_1$ in $Z_3$ can be contracted to non-singular rational curves 
simultaneously.
For this we first consider the divisor $E_1$  in $Z_1$, so that $E_1\simeq\mathbb{CP}^1\times\mathbb{CP}^1$, and take the cohomology class of a fiber of the projection to $\mathbb{CP}^1$ which is different from the projection to $C_1$.
Next pullback the class by the blowup $\mu_2$ and push it to $E_1\subset Z_3$ by $\mu_3$.
(In Figure \ref{fig:trans} these cohomology classes are represented by non-dotted lines   on $E_1$.)
Thus we obtain a cohomology class on $E_1\subset Z_3$ whose self-intersection number is zero.
The linear system on this $E_1$ having this cohomology class is clearly a free pencil, and induces a morphism to $\mathbb{CP}^1$.
Let $g:E_1\to\mathbb{CP}^1$ be this morphism.
General fibers of $g$ are non-singular rational curves, and there exist precisely 2 singular fibers, both of which consist of 2 non-singular rational curves intersecting at a point.
The same is true for $\ol E_1$, and let $\ol g:\ol E_1\to \mathbb{CP}^1$ the morphism corresponding to $g$.
Now $g$ and $\ol g$ naturally fit on
the intersection $E_1\cap \ol E_1$ and form a morphism $g\cup \ol g:E_1\cup \ol E_1 \to \mathbb{CP}^1\cup\mathbb{CP}^1$.
Here note that these two $\mathbb{CP}^1$-s are identified at 2 points,
and $g\cup\ol g$ has reducible fibers exactly over these 2 points, both of which consist of {\em three}
rational  curves.
In Figure \ref{fig:trans} (d), these 2 reducible fibers are written by 3 bold lines respectively.

We are going to show that the reducible  connected divisor $E_1\cup\ol E_1$ on $Z_3$ can be contracted along $g\cup\ol g$.
For this we need to examine the normal bundle, $[E_1+ \ol E_1]|_{E_1\cup\ol E_1}$.
The restriction of the normal bundle $N_{E_1/Z_3}=[E_1]|_{E_1}$ is degree $(-2)$ on irreducible fibers of $g$, and degree $(-1)$ on the 2 irreducible components of the (two) singular fibers. (See Figure \ref{fig:trans} (d).)
From this we deduce that the restriction of the line bundle $[E_1+ \ol E_1]|_{E_1\cup\ol E_1}$ is $(-2)$ on irreducible fibers of $g\cup\ol g$, and $(-1)$ on the end components of the reducible fibers, while $(-1)+(-1)=-2$ on the middle component
of the reducible fibers.
Now by the relative version of Nakai-Moishezon criterion for ampleness,
these numerical data imply that the dual line bundle $[E_1+ \ol E_1]^*|_{E_1\cup\ol E_1}$ is  $(g\cup \ol g)$-ample.
Moreover again from the  numerical data,  for a direct image of the dual bundle, we have 
$$
R^1(g\cup\ol g)_*([E_1+ \ol E_1] ^* | _{E_1\cup\ol E_1} ) ^ {\otimes m}=0 \hsp{\text{for any}}\hsp m>0.
$$
Therefore by a theorem of Fujiki \cite[Theorem 2]{F75}, the divisor $E_1\cup \ol E_1$ can be contracted to $\mathbb{CP}^1\cup\mathbb{CP}^1$ along the morphism $g\cup\ol g$.
Let $\mu_4:Z_3\to Z_4$ be the birational morphism obtained this way, and put $l_4:=\mu_4(E_1)$ and $\ol l_4:=\mu_4(\ol E_1)$, so that $l_4\cup\ol l_4$ can be naturally identified 
with the target space of $g\cup\ol g$.
As the degree of the restriction of $N_{E_1/Z_3}$ to irreducible fibers of $g$ is $(-2)$ as above, we have
$
{\rm{Sing}}\, Z_4 = l_4\cup\ol l_4,
$
and possibly outside the 2 points $l_4\cap \ol l_4$,
$Z_4$ has ordinary double points along $l_4\cup\ol l_4$.
(In Section \ref{ss:blowup} we will obtain an explicit defining equation of $Z_4$ around the 2 points.)
This way by contracting $E_1\cup\ol E_1$ in $Z_3$ we have obtained a singular variety $Z_4$.
Then the morphism $\Phi_2$ descends to $Z_4$:

\begin{proposition}\label{prop:factor1}
Let $\Phi_2:Z_2\to Y$ be the generically 2 to 1  covering as before. 
Then $\Phi_2$ descends to a morphism $Z_4\to Y$.
Namely there is a morphism $\Phi_4:Z_4\to Y$ such that $\Phi_4\circ\mu_4\circ\mu_3=\Phi_2$.
\end{proposition}

\proof
We first show that $\Phi_2$ descends to a morphism $\Phi_3:Z_3\to Y$.
Recall that $\Phi_2$ is induced by the system $|\mathscr L_2|$, where $\mathscr L_2=\mu_2^*\mathscr L_1$ and $\mathscr L_1=\mu_1^*2F-E_1-\ol E_1$.
In accordance with those on $\Sigma_1$ in the proof of Proposition \ref{branchdivisor1},
let $(0,1)$ be the fiber class of the projection $E_1\to C_1$.
Then we obtain 
$
\mathscr L_1|_{E_1}\simeq 2\mu_1^*(F|_{C_1})-N_{E_1/Z_1}
\simeq 2(\mu_1^*K_S^{-1}|_{C_1})-\mathscr O(-1,-2)
\simeq 2\mu_1^*\mathscr O_{C_1}(-1)+\mathscr O(1,2)
\simeq\mathscr O(1,0).
$
(See Figure \ref{fig:trans} (b) for $N_{E_1/Z_1}\simeq
\mathscr O(-1,-2)$.)
Further, as the curve $C_2\subset Z_1$ intersects $E_1$  transversally at exactly 1 point and the same for $\ol E_1$
(again see Figure \ref{fig:trans} (b)), we have
$
\mathscr L_1|_{C_2}\simeq (\mu_1^*2F)|_{C_2} -\mathscr O_{C_2}(2)
\simeq 2F|_{C_2} - \mathscr O_{C_2}(2)
\simeq 2K_S^{-1}|_{C_2} - \mathscr O_{C_2}(2)
\simeq \mathscr O_{C_2}.
$
Therefore pulling back to  $Z_2$ we obtain
$
\mathscr L_2|_{E_1} \simeq \mu_2^*\mathscr O(1,0),
\quad
\mathscr L_2|_{E_2} \simeq \mathscr O_{E_2},
$
and analogous result for the restrictions to $\ol E_1$ and $\ol E_2$.
These imply that the direct image sheaf $(\mu_3)_*\mathscr L_2=:\mathscr L_3$ is still invertible and $\mathscr L_3|_{E_1}\simeq \mu^*_2\mathscr O_{E_1}(1,0)$.
If we use the projection $g:E_1\to l_4$, the last isomorphism can be rewritten as
\begin{align}\label{L3}
\mathscr L_3|_{E_1}\simeq g^*\mathscr O_{l_4}(1).
\end{align}
Then since $\mathscr L_2\simeq \mu_3^*\mathscr L_3$,
the morphism associated to $|\mathscr L_2|$ factors through the morphism associated to $|\mathscr L_3|$.
Letting $\Phi_3$ be the last morphism, 
this means $\Phi_2=\Phi_3\circ\mu_3$ as claimed.

In a similar way we next show that $\Phi_3$ descends to a morphism $\Phi_4:Z_4\to Y$.
 From \eqref{L3} the direct image $(\mu_4)_*\mathscr L_3=:\mathscr L_4$ is still an invertible sheaf (on $Z_4$),
whose restriction to $l_4$ is of degree 1.
 Then by the natural isomorphisms $H^0(Z_2,\mathscr L_2)\simeq H^0(Z_3,\mathscr L_3)\simeq 
 H^0(Z_4,\mathscr L_4)$ the map associated to $|\mathscr L_3|$ factors through the map induced by 
 $|\mathscr L_4|$.
 Therefore if we define $\Phi_4$ to be the map associated to $|\mathscr L_4|$, 
 we have $\Phi_3=\Phi_4\circ\mu_4$, as claimed.
 Hence we have obtained $\Phi_2=\Phi_3\circ\mu_3=\Phi_4\circ\mu_4\circ\mu_3$.
 \proofend

\vspace{3mm}
Thus we arrived at the following situation:
\begin{equation}\label{cd3}
\xymatrix{
   Z_2 \ar@{->}[d]_{\mu_2}   \ar@{->}[r]^{\mu_4\circ\mu_3}  \ar@{->}[rd]^{\Phi_2} & 
   Z_4 \ar@{->}[d]^{\Phi_4}  \\
 Z_1 \ar@{->}[r]_{\Phi_1} & Y 
}
 \end{equation}
 where all maps are morphisms and the 2 triangles are commutative.
Since $\mu_4\circ\mu_3$ is birational, $\Phi_4$ is still a degree 2 morphism branching at the divisor $B$.
But in contrast with $\Phi_2$, it does not contract divisors anymore:

\begin{proposition}\label{prop:phi4}
The morphism $\Phi_4$ does not contract any divisor to a point or a curve.
\end{proposition}

\proof
It is enough to show that the morphism $\Phi_1:Z_1\to Y$ does not contract any irreducible divisor
 other than $E_1$ and $\ol E_1$.
Let $D$ be such a divisor.
If $D$ is real,
then $D\in |\mu_1^*(kF)-lE_1-l\ol E_1|$ for some $k\ge 1$ and $l\ge 0$,
and $(\mu_1^*2F-E_1-\ol E_1)^2\cdot D =0$ by the contractedness property.
For computing this intersection number, 
we notice, as $\mu_1^*F|_{E_1} \simeq \mu_1^* (K_S^{-1}|_{C_1}) \simeq \mathscr O_{E_1}(0,-1)$, that we have
$(\mu_1^*F)^2 \cdot E_1 = (\mu_1^*F|_{E_1})^2 = 0$.
Hence we have
\begin{align*}
(\mu_1^*2F-E_1-\ol E_1)^2 \cdot E_1 & =  4 \mu_1^* F^2 \cdot E_1 - 4 \mu_1^* F \cdot (E_1+\ol E_1) \cdot E_1
+ (E_1+\ol E_1)^2\cdot E_1\\
&= 4 \cdot 0 - 4 \mu_1^*F \cdot E_1^2 + E_1^3 \\
&= -4 \mathscr O_{E_1}(0,-1) \cdot \mathscr O_{E_1} (-1,-2) + \mathscr O_{E_1}(-1,-2)^2 \\
&= -4+4=0,
\end{align*}
and the same for $\ol E_1$.
From these, recalling $F^3=0$ (as we are over $4\mathbb{CP}^2$), we obtain
\begin{align*}
(\mu_1^*2F-E_1-\ol E_1)^2\cdot D &= (\mu_1^*2F-E_1-\ol E_1)^2 \cdot (\mu_1^*(kF)-lE_1-l\ol E_1) \\
&= (\mu_1^*2F-E_1-\ol E_1)^2 \cdot \mu_1^*(kF) \\
&= E_1^2\cdot \mu_1^*(kF) + \ol E_1^2\cdot \mu_1^*(kF) - 4k \mu_1^*F^2\cdot (E_1+\ol E_1) \\
&= 2 N_{E_1/Z_1} \cdot \mu_1^*(kF)|_{E_1} \\
&= 2k \mathscr O_{E_1}(-1,-2)\cdot \mathscr O_{E_1}(0,-1) = 2k.
\end{align*}
Therefore we have $(\mu_1^*2F-E_1-\ol E_1)^2\cdot D >0$.
Hence $D$ is not contracted to a curve or a point by $\Phi_1$.
When $D$ is not real, by applying the above computations
for $D+\ol D$ instead of $D$,
we obtain $(\mu_1^*2F-E_1-\ol E_1)^2\cdot (D+ \ol D) >0$.
Hence, since $(\mu_1^*2F-E_1-\ol E_1)^2\cdot D
= (\mu_1^*2F-E_1-\ol E_1)^2\cdot \ol D$, 
we again conclude $(\mu_1^*2F-E_1-\ol E_1)^2\cdot D >0$.
Hence in the non-real case too, $D$ cannot be contracted to a point or a curve by $\Phi_1$, as claimed.
\proofend

\vsp
As a consequence, we obtain the following
\begin{proposition}\label{prop:isolated}
The branch divisor $B$ has only isolated singularities.
\end{proposition}

\proof
Let $Z_4\stackrel{\mu_5}{\to}Z_5\stackrel{\Phi_5}\to Y$ be the Stein factorization of the morphism $\Phi_4$. 
$\mu_5$ is necessarily birational.
Then $\Phi_5$ is just a double covering with branch
$B$. Hence if $B$ has singularities along a curve, so is $Z_5$.
Since $B$ does not contain $l$ and ${\rm Sing}\, Z_4=l_4\cup\ol l_4$, this means that 
the birational morphism $\mu_5$ resolves the singularities along the curve. Hence
$\mu_5$ contracts a divisor. 
This contradicts Proposition \ref{prop:phi4}.
\proofend
\vsp

We note that the proof of Proposition \ref{prop:phi4} means that the original anticanonical map 
$\Phi:Z\to Y$ does not contract any divisor.
We also note that the morphism $\mu_5$ in the
proof of 
Proposition \ref{prop:isolated} contracts
(at worst) finitely many curves, and
all these curves are over  singular points of  $B$.
We will investigate these singularities 
in detail in Section \ref{ss:sing}.

\section{Defining equation of the branch quartic hypersurface}
\label{s:defeq}
In the last section we analyzed the anticanonical system on the twistor space in detail and obtained
the space $Z_4$ and a degree 2 morphism $\Phi_4:Z_4\to Y$ which is explicitly birational  to the original anticanonical map, and which does not contract any divisor.
We further showed that the branch divisor $B$ of $\Phi_4$ is a cut of $Y$ by a quartic hypersurface in $\mathbb{CP}^4$.
In this section we shall determine defining equation of this quartic hypersurface.
We also determine the number of singularities of the branch divisor.

\subsection{Finding double curves on $B$}\label{ss:dc}
Our way for obtaining the equation includes finding hyperplanes $H\subset \mathbb{CP}^4$ such that, regarding 
the intersection $H\cap Y$ 
(which is either a plane or a cone as in the beginning of 
Section \ref{ss:modify}) as a reduced divisor on $Y$, the restriction $B|_{H\cap Y}$ is a double curve (i.e.\,a non-reduced curve of multiplicity 2).
So it is similar to the method of Poon \cite{P92} (for the case of $3\mathbb{CP}^2$),
but the origin of some of the double curves is different from the case of $3\mathbb{CP}^2$.

We keep the  notations from the last section.
We are going to show the existence of {\em five} double curves, two of which are easy to find as we see now.
First we recall there are diagrams in \eqref{cd1};
$\Lambda$ is a conic in $\mathbb{CP}^2$ and for each $\lambda\in \Lambda$, $S_{\lambda}:=f^{-1}(\lambda)$ is a member of the pencil $|F|$, i.e.\,$\Lambda$ is a parameter space of $|F|$.
For any $\lambda\in\Lambda$,  $\pi^{-1}(\lambda)$
is a plane containing the line $l$.
Let $S_1=S_1^++S_1^-$ and $S_2= S_2^++S_2^-$ be the reducible members as in \eqref{reduciblemembers}, and let $0$ and $\infty$ be the points of $\Lambda$
such that $f^{-1}(0)=S_1$ and $f^{-1}(\infty)=S_2$ hold.
Let $ H_1$ and $ H_2$ be the hyperplanes in $\mathbb{CP}^4$ which are the inverse images of the tangent line of $\Lambda$ at the points $0$ and $\infty$ respectively under the projection $\pi$.
Then the restrictions $ H_1|_Y$ and $ H_2|_Y$ are double planes, 
and we have $\Phi^{-1}( H_1)=2S_1$ and $\Phi^{-1}( H_2)=2S_{2}$.
Letting $L_1=S_1^+ \cap S_1^-$ and $L_{2}=S_{2}^+\cap S_{2}^-$ be the twistor lines as before, we  define 
\begin{align}\label{dc}
\mathscr C_1:=\Phi(L_1) \,\text{ and }\, \mathscr C_2:=\Phi(L_{2}).
\end{align}
Then since  $\Phi$ is a real map and $0$ and $\infty$ are real points,
and since $\Phi$ does not contract any divisor by Proposition \ref{prop:phi4},
$S_1^+$ and $S_1^-$ are mapped birationally to  the plane $\pi^{-1}(0)$ and 
$S_2^+$ and $S_2^-$ are mapped birationally to the plane $\pi^{-1}(\infty)$.
Hence $\mathscr C_1$ and $\mathscr C_2$ are contained in the branch divisor $B$ in such a way that the restrictions $B|_{ H_1\cap Y}$ and $B|_{ H_2 \cap Y}$ respectively contain $\mathscr C_1$ and $\mathscr C_2$ by multiplicity 2.
But as we know that $B$ is a cut of $Y$ by a hyperquartic surface by Proposition \ref{branchdivisor1}, $\mathscr C_1$ and $\mathscr C_2$ must be conics.
So we call these two {\em double conics}.

These double conics are analogous to the 4 conics contained in the coordinates tetrahedron in $\mathbb{CP}^3$ used in \cite{P92} in the case of $3\mathbb{CP}^2$, but in the present case  there are only 2 since there are only 2 reducible members of $|F|$.
Next we find other 3 double curves.
For this recall that our twistor space $Z$ contains the surface $S$ constructed in Section \ref{s:1} as a real member of $|F|$.
Take up the blowup $\epsilon:S\to \mathbb{CP}^1\times\mathbb{CP}^1$ which was concretely given when constructing $S$,
and let $e_i$ and $\ol e_i$, $1\le i\le 4$, be the exceptional curves of the blowup.
Here we take indices such that $e_i\cdot C_1=1$ for $i=1,2,3$ 
and $e_4\cdot C_2=1$.
Next let $\{\alpha_i\set 1\le i\le 4\}$ be an orthonormal basis of $H^2(4\mathbb{CP}^2,\mathbb Z)$ determined by $\{e_i\}$; namely
letting $t:Z\to 4\mathbb{CP}^2$ be the twistor fibration, $t^*\alpha_i|_S=e_i-\ol e_i$ in $H^2(S,\mathbb Z)$.
Then we have the following proposition which is 
technically significant for our purpose:

\begin{proposition}\label{prop:redivisor}
For each $i$ with $1\le i\le 3$ (not 4), 
the linear systems $|F+t^*\alpha_i|$ and
$|F-t^*\alpha_i|$ consist of a single element.
Moreover, all these 6 divisors are irreducible.
\end{proposition}

\proof
In this proof for simplicity we write $\alpha_i$ for $t^*\alpha_i$.
It is enough to prove the claim for the system $|F+\alpha_i|$.
Fix any $i$ with $1\le i\le 3$.
We have $(F +   \alpha_i ) | _S = K_S^{-1} + ( e_i - \ol e_i ) =
\epsilon^* \mathscr O(2,2) - \sum _ {j=1} ^4 ( e_j + \ol e_j ) +  ( e_i - \ol e_i )$, which can be rewritten as
$$
 \Big( \epsilon^* \mathscr O(1,0) - \sum_{1\le j\le3} \ol e_j \Big) +
 \Big( \epsilon^* \mathscr O(1,2)   - \ol e_i - e_4 - \ol e_4 - \sum _ {1\le j\le 3,\, j\neq i }  e_j \Big).
$$
The intersection number between 
$\ol C_1\sim \epsilon^* \mathscr O(1,0) - \sum _ {j=1}^3 \ol e_j$ and
the above $(F+\alpha_i)|_S$ is easily 
 computed to be $(-2)$.
Hence $\ol C_1$ is a fixed component of $|(F +   \alpha_i ) | _S|$.
Further counting dimension, the remaining system $| \epsilon^* \mathscr O(1,2)   - \ol e_i - e_4 - \ol e_4 - \sum _ {j=1,\, j\neq i } ^ 3 e_j | $ consists of a single member, which is the strict transform of a $(1,2)$-curve passing through the 5 points $\ol q_i, q_4, \ol q_4$ and $q_j$ with $1\le j\le 3$ and $j\neq i$,
where $q_i=\epsilon(e_i)$ and $\ol q_i = \epsilon (\ol e_i)$.
Thus we have $h^0 ( ( F +   \alpha_i ) | _S )  = 1$.
On the other hand from Riemann-Roch formula and Hitchin's vanishing theorem \cite{Hi80} we deduce $H^1(Z, \alpha_i)=0$.
Then by the standard exact sequence $ 0 \to   \alpha_i  \to F +  \alpha_i  \to ( F +   \alpha_i ) |_S
\to 0$, we obtain the exact sequence  
$ 0 \to H^0(Z,  \alpha_i )  \to H^0 (Z, F  +  \alpha_i )  \to H^0 (S, ( F +   \alpha_i ) |_S)
\to 0$.
As $H^0(Z,  \alpha_i)=0$, we obtain $ H^0 (Z, F  +  \alpha_i )  \simeq H^0 (S, ( F +   \alpha_i ) |_S)$.
Hence we get $ H^0 (Z, F  +  \alpha_i )  \simeq \mathbb C$.

For the irreducibility, all the divisors $S_1^+, S_1^-, S_{2}^+$ and $S_{2}^-$ 
(the irreducible components of reducible members of $|F|$) are degree 1  on $Z$.
Moreover, it is not difficult to show that the Chern classes of these divisors are given by {\em the half of} the following classes:
\begin{align}\label{cc1}
F -  \sum _{1\le j\le 4}   \alpha_j, \hsp
F + \sum _{1\le j\le 4}   \alpha_j, \hsp
F + \alpha _4 - \sum _{1\le j\le 3}  \alpha_j  , \hsp
F - \alpha _4 + \sum _{1\le j\le 3}  \alpha_j .
\end{align}
Then it is a easy to check that a sum of any two  of these 4 classes (allowing
to choose the same one) are not equal to 
$2(F+ \alpha_i)$, for any $1\le i\le 3$.
This implies the desired irreducibility.
(On the other hand, by \eqref{cc1}, the systems $|F\pm\alpha_4|$ are also non-empty, but both of them consist of a single {\em reducible} member.)
\proofend

\vspace{3mm}
In the following for $1\le i\le 3$ we denote by $X_i$ for the unique member of $|F+ t^*\alpha_i|$.
Then $\ol X_i\in |F-t^*\alpha_i|$, and $X_i + \ol X_i\in |2F|$.
Thus we obtained 3 reducible real members of the anticanonical system on $Z$.
We remark that from the proof of Proposition \ref{prop:redivisor} these 3 members originally come from the choice of 3 points on $C_1$ in the construction of $S$ at the beginning of Section \ref{s:1}.
By using these we obtain a special basis of $H^0(2F)\simeq\mathbb C^5$ as follows:

\begin{proposition}\label{prop:basis4}
For any $i$ with $1\le i\le 3$, let $\xi_i\in H^0(Z,2F)$ be an element such that $(\xi_i)= X_i + \ol X_i$.
Then $S^2H^0(Z,F)$ $(\simeq\mathbb C^3)$ and any two among $\{\xi_i
\set 1\le i\le 3\}$ generate  $H^0(2F)\,(\simeq\mathbb C^5)$.
 \end{proposition}

\proof
Let $S\in |F|$ be any real irreducible member and take  $s_0\in H^0(F)$ with $(s_0)=S$.
Let $s_1\in H^0(F)$ be any element satisfying $s_1\not\in\mathbb Cs_0$.
Then $\{s_0,s_1\}$ is a basis of $H^0(F)$ and $\{s_0^2,s_0s_1,s_1^2\}$ is a basis of $S^2H^0(F)$.
We consider the  exact sequence
\begin{align}\label{ses:basic5}
0 \lra H^0 (F) \stackrel{\otimes s_0}{\lra} H^0 (2F) \lra H^0(2K_S\inv) \lra 0
\end{align}
appeared in the proof of Proposition \ref{prop:2F}.
For proving the claim of the proposition, it suffices to show that for any subset $\{i,j\}\subset\{1,2,3\}$,
the images of $s_1^2, \xi_i, \xi_j$ by the restriction map to $S$ generate $H^0(2K_S^{-1})$.
For this, the divisor $(s_1^2|_S)$ is exactly $2C$, where $C$ is the unique anticanonical curve (i.e.\,the cycle of 4 rational curves).
On the other hand we have $(\xi _i | _S)=X_i|_S + \ol X_i |_S$, and from the proof of Proposition \ref{prop:redivisor} we know the curves
$X_i|_S$ and $ \ol X_i |_S$ in concrete forms,
and it is not difficult to verify that the 3 bi-anticanonical curves $2C, X_i|_S+ \ol X_i |_S, X_j|_S+ \ol X_j |_S$ are linearly independent.
 Hence the 3 images generate $H^0(2K_S^{-1})$.
 \proofend
 
 \vspace{3mm}
In the sequel for obtaining nice coordinates, we choose a slightly different basis $\{u_1,u_2\}$ of $H^0(Z,F)$ as follows.
Namely respecting the reducible members,  we choose those satisfying $(u_1) = S_1$ and $(u_{2})= S_2$.
By Proposition \ref{prop:basis4} the collection $\{u_1u_{2}, u_1^2, u_{2}^2, \xi_1, \xi_2\}$ is a basis of $H^0(Z,2F)$.
The target space of the anticanonical map $\Phi:Z\to \mathbb{CP}^4$ is nothing but 
the dual projective space $\mathbb P H^0(Z,2F)^*\,(\simeq\mathbb{CP}^4)$,
and if we put
\begin{align}\label{z_i}
z_0: = u_1u_{2}, \,\, z_1 := u_1^2, \,\, z_2:=u_{2}^2,\,\,
z_3 := \xi _1, \,\, z_4 := \xi _2,
\end{align}
then $(z_0,z_1,z_2,z_3,z_4)$ can be used as homogeneous coordinates on it.
(Here we remark that there is no special reason to choose $\xi_1$ and $\xi_2$.
Any choice of two among $\{\xi_1,\xi_2,\xi_3\}$ leads to the same description below.)
As $\mathbb{CP}^4=\mathbb PH^0(2F)^*$ as above, $\mathbb{CP}^4$
is equipped with a real structure and by \eqref{z_i} it is
just the complex conjugation with respect to the above coordinates.
In these coordinates the  scroll $Y=\Phi(Z)$ is explicitly defined by the equation
\begin{align}\label{eqn:Y1}
z_0^2 = z_1 z_2,
\end{align}
and the ridge $l$ of $Y$ is given by $z_0=z_1=z_2=0.$
For each $0\le i\le 4$ we define a hyperplane  by $H_i:=\{z_i=0\}$.
Obviously $l\subset H_i$ for $i=0,1,2$, and $l\not\subset H_i$ for $i=3,4$.
In particular $H_1|_Y$ and $H_2|_Y$ are double planes, $H_0|_Y$ is the sum of these 2 planes,
and $H_3| _ Y$ and $H_4 |_ Y$ are cones whose vertices are the  points $H_3\cap l$ and $H_4\cap l$ respectively.

Let $z_5\in H^0(Z,2F)$ be an element such that $(z_5)= X_3 + \ol X_3$.
Then by Proposition \ref{prop:basis4} we can write
\begin{align}\label{z5}
z_5= a_0z_0 + a_1 z_1 + a_2 z_2 + a_3z_3 + a_4 z_4
\end{align}
for some $a_i\in\mathbb R$.
Let $H_5:=\{z_5=0\}$. 
Since  $z_5\not\in S^2H^0(Z,F)$ clearly, $H_5|_Y$ is also a cone.

Now the following proposition provides the promised 3 double curves on $B$:
\begin{proposition}\label{prop:dq1}
For $i=3,4,5$, the intersection of the branch divisor $B$ with the cone  $H_i\cap Y$ is a double curve of $B$. 
\end{proposition}

\proof
Let $i$ be any one of $3,4,5$.
By definition of the hyperplane $H_i$, we have $\Phi^{-1}(H_i)= X_{i-2} + \ol X_{i-2}$.
Since $\Phi$ does not contract any divisor by Proposition \ref{prop:phi4}, this implies $\Phi(X_{i-2})=\Phi( \ol X_{i-2})=H_i\cap Y$.
As $\Phi$ is degree 2, this means that the restrictions $\Phi |_{X_{i-2}}$ and $\Phi |_{\ol X_{i-2}}$ are birational over the cone $H_{i}|_Y$.
(In particular, the non-real degree 2 divisors $X_{i-2}$ and $\ol X_{i-2}$ are birational to a cone.)
Therefore the curve $X_{i-2}\cap \ol X_{i-2}$ is the ramification divisor of the restriction
$\Phi|_{\Phi^{-1}(H_i)}:\Phi^{-1}(H_i)  \to H_i\cap Y$,
and hence 
$$
\mathscr C_3 := \Phi ( X_1 \cap \ol X_1) \quad  \mathscr C_4:= \Phi ( X_2 \cap \ol X_2)\hsp{\text{and}}\hsp
\mathscr C_5 := \Phi ( X_3 \cap \ol X_3)
$$
are branch divisors  when restricted to $\Phi^{-1}(H_i)$. 
This implies the claim of the proposition.
\proofend
\vspace{3mm}

As in the proof, we use the letters $\mathscr C_3, \mathscr C_4$ and $\mathscr C_5$ to mean the 3 double curves in the proposition.
Then because we know that $B$ is a cut of $Y$ by a quartic hypersurface, we have $\mathscr C_i\in |\mathscr O_{Y\cap H_{i}}(2)|$ , where $\mathscr O_{Y\cap H_i}(2)
:= \mathscr O_{\mathbb{CP}^4}(2)|_{Y\cap H_i}$.
Namely, $\mathscr C_3, \mathscr C_4$ and $\mathscr C_5$
are intersection of the cone $Y\cap H_i$ with a
quadratic in $H_i=\mathbb{CP}^3$.
From this it follows that these 3 curves are of degree 4 in $\mathbb{CP}^4$.
So in the following we call these double curves  {\em double quartic curves}.
From our choice of the coordinates, any double curves can be written as , as sets,
\begin{align}\label{dc3}
\mathscr C_i = B \cap H_i, \quad 1\le i\le 5.
\end{align}
Since $B$ and $H_i$ are real, all $\mathscr C_i$-s are real curves.

\subsection {Quadratic hypersurfaces containing the double curves}\label{ss:quadric}
In this section we show that there exists a quadratic hypersurface in $\mathbb{CP}^4$ which contains the double conics $\mathscr C_1$, $\mathscr C_2$ and the double quartic curves $\mathscr C_3, \mathscr C_4$ and $\mathscr C_5$, and also
show that such hyperquadric is unique up to the defining equation of the scroll $Y$.

\begin{figure}
\includegraphics{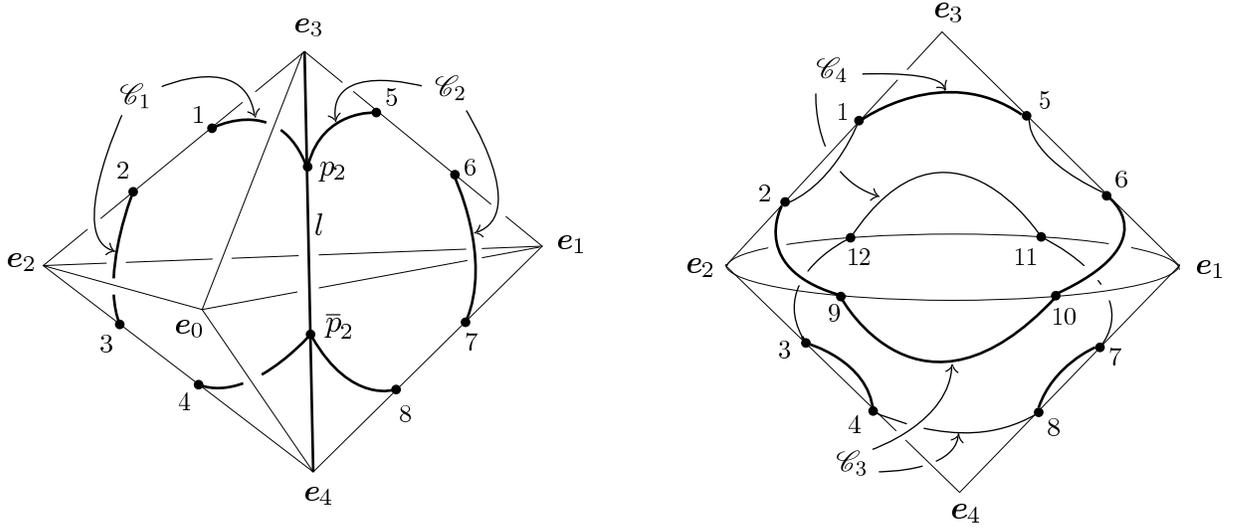}
\caption{The double curves of the branch divisor $B$, except $\mathscr C_5$. 
Both pictures lie on the same $\mathbb{CP}^4$,
and the common numbers represent the same point, all of which
are ordinary double points of $B$.
The left picture indicates all coordinate hyperplanes, 2-planes and lines, as well as double conics.
In the right picture the upper and lower halves are the cones $Y\cap H_4$ and $Y\cap H_3$ respectively, on which 
the double quartic curves
$\mathscr C_4$ and $\mathscr C_3$ lie.}
\label{fig:dc}
\end{figure}

First we make it clear how the 5 double curves of $B$ intersect each other.
For this for each $0\le i\le 4$ we denote by $\bm e_i\in\mathbb{CP}^4$ for the point whose coordinates are zero except $z_i$-component,
and
define some lines  as follows:
for each pair $(i,j)$ with $i=1,2$ and $j=3,4,5$, define $l_{ij} := H_i\cap Y\cap H_j  $.
Since $H_i\cap Y$ is a (double) plane for $i=1,2$, this is a line.
Thus we get 6 lines.
If $j\neq 5$, these are coordinate lines and 
in Figure \ref{fig:dc}, $l_{14}= \ol{\bm e_2 \bm e_3},\,
l_{13}=\ol{\bm e_2 \bm e_4},\,l_{24}=\ol{\bm e_1 \bm e_3}, \,
l_{23} = \ol{\bm e_1 \bm e_4}$.
(We do not write pictures of $l_{15}$ and $l_{25}$ because these are not
coordinate lines. But this is just a matter of a choice of coordinates and these two play the same role as other 4 lines.)
Also, for the ridge $l$ we have $l=\ol{\bm e_3\bm e_4}$
(the bold line on the left picture in Figure \ref{fig:dc}).

Then since the intersection of the branch divisor $B$ with the plane $H_i\cap Y=\{z_0=z_i=0\}$ $(i=1,2)$ is the double conics $\mathscr C_i$ and since the line $l_{ij}$ ($3\le j\le 5$) is contained in this plane,  the intersection $B\cap l_{ij}$ consists of (not 4 but) 2 points,
and $B\cap l_{ij}= \mathscr C_i\cap l_{ij}$.
Moreover, as $\mathscr C_i=\Phi(L_i)$, these 2 points cannot be identical.
Thus for each of the 6 lines $l_{ij}$,  $B\cap l_{ij}$ consists of  2 points.
(In Figure \ref{fig:dc} these points are represented by
numbered points $1,2,\cdots,7,8$.)
On the other hand, we have $B\cap H_j=\mathscr C_j$.
Hence as $l_{ij}\subset H_j$, we obtain $B\cap l_{ij}\subset
\mathscr C_j\cap l_{ij}$ for $i=1,2$ and $j=3,4,5$.
But since $\mathscr C_j$ is an intersection of $Y\cap H_j$ with 
a quadric surface in $H_j$, $\mathscr C_i\cap l_{ij}$ consists
at most 2 points.
Hence we have the coincidence $B\cap l_{ij}=
\mathscr C_j\cap l_{ij}$ for these $ i$ and $j$.
Therefore
we have $B\cap l_{ij} = \mathscr C_i\cap \mathscr C_j$
for $i=1,2$ and $j=3,4,5$.
By a similar reason, the intersection $B\cap l$ also consists of 2 points, which are exactly $\mathscr C_1\cap\mathscr C_2$.
In Figure \ref{fig:dc} these points are displayed as $p_2$ and $\ol p_2$.
On the other hand, 
 for each pair $(j,k)$ with $3\le j<k\le 5$  we define 
a plane  $P_{jk}$ by $P_{jk}= H_j\cap H_k$.
Then since $B\cap H_j$ is contained in a quadric surface, and $Y\cap P_{jk}$ is a conic,
$B\cap P_{jk}$ ($3\le j<k\le 5$) consists of 4 points, and it coincides with $\mathscr C_j\cap \mathscr C_k$.
In Figure \ref{fig:dc}, for the case $(j,k)= (3,4)$, these are  represented by numbered points $9,10,11,12$.
(For avoiding confusion we do not write a picture for $\mathscr C_3\cap\mathscr C_5$ and $\mathscr C_4\cap\mathscr C_5$.
 The way how these curves intersect is completely analogous to
that of $\mathscr C_3$ and $\mathscr C_4$.)
We list all these intersections: 
\begin{itemize}
\item 2 points $\mathscr C_1\cap\mathscr C_2$, which are exactly $p_2$ and $\ol p_2$,
\item $12$ points $\mathscr C_i\cap\mathscr C_j$ with $i=1,2$ and $j=3,4,5$,
\item $12$ points $\mathscr C_3\cap\mathscr C_4$,
 $\mathscr C_3\cap\mathscr C_5$ and
 $\mathscr C_4\cap\mathscr C_5$.
\end{itemize}
Collecting these, we obtain 26 points in total.
Since all the double curves are the image of curves in $Z$ by a map
which is degree 1 on these curves,
these 26 points form 13 conjugate pairs.
Among these 26 points, the 2 points $\mathscr C_1\cap \mathscr C_2$ are on the singular locus $l$ of $Y$,
and other 24 points are ordinary double points of $B$.
(In some sense these 24 points are analogous to the 12 ordinary double points of the branch quartic surface appeared in \cite{P92} and \cite{KK92}.) In Section \ref{ss:sing}, we will show that the 2 points
$p_2$ and $\ol p_2$ are $A_3$-singular points of $B$.

With these situation in hand,
 we next show the existence of a hyperquadric which contains all the double conics:

\begin{proposition}\label{prop:Q}
There exists a real quadratic hypersurface in $\mathbb{CP}^4$ which contains all the 5 double curves $\mathscr C_i$, $1\le i\le 5$, and which is different from the scroll $Y$.
Moreover, such a hyperquadric is unique in the following sense:
if $Q$ and $Q'$ are defining quadratic polynomials of two such hyperquadrics,
then there exists $(c,c')\in\mathbb R^2$ with $(c,c')\neq (0,0)$ such that $cQ-c'Q'\in ( z_0^2-z_1z_2)$.
(Note that since the scroll \, $Y$ contains all 
the double curves, presence of this ambiguity is obvious from the beginning.)
\end{proposition}

\proof As we have already seen,
for $i=3,4$ the intersection $H_i\cap Y$ is a quadratic cone in $H_i=\mathbb{CP}^3$, and the double quartic curve $\mathscr C_i$ belongs to $|\mathscr O_{H_i\cap Y}(2)|$.
In the above homogeneous coordinates the intersection $H_3\cap H_4$ is a plane defined by $z_3=z_4=0$, 
and $H_3\cap H_4\cap Y$ is a conic defined by $z_0^2=z_1z_2$, on which the 4 points $\mathscr C_3\cap \mathscr C_4$ lie.
Conics on the plane passing through these 4 points form a pencil, which 
is invariant under the real structure.
Choose any real one of such conics, and let $q(z_0,z_1,z_2)$ be its defining equation with real coefficients, which is of course uniquely determined up to rescaling.
Among the above 26 points  there are exactly 8 points lying on $H_3$
(which are the points $3,4$ and $7$ to $12$ in Figure \ref{fig:dc}), four of which are the above 4 points
on $H_3\cap H_4$ (the points $9,10,11,12$ in Figure \ref{fig:dc}).
Any quadratic polynomial  on $H_3$ whose restriction to $H_3\cap H_4$ equals $q$ is of the form
\begin{align}
Q_3 = q(z_0,z_1,z_2) + a_0z_0z_3 + a_1z_1z_3 + a_2z_2z_3
+ a_3z_3^2.
\end{align}
Imposing that the quadric $(Q_3)$ passes the remaining 
4 points ($3,4,7,8$), the coefficients $a_0,a_1,a_2$ and $a_3$ are uniquely determined (without an ambiguity of rescaling), and 
they are real since the set of 4 points (3,4,7,8) is real.
Then since elements of $|\mathscr O_{H_3\cap Y}(2)|$ which 
go through the 8 points is unique by dimension counting, it follows that the quadratic surface $(Q_3)$ automatically contains $\mathscr C_3$.
The situation is the same for $H_4$, and let 
\begin{align}
Q_4 = q(z_0,z_1,z_2) + b_0z_0z_4 + b_1z_1z_4 + b_2z_2z_4
+ b_3z_4^2.
\end{align}
be the  quadratic polynomial on $H_4$ with  real coefficients,
which is uniquely determined 
by $q$ and the 8 points ($ 1,2,5,6$ and $9$ to $12$) on $H_4$.
Then  $(Q_4)\supset \mathscr C_4$.

Since $Q_3|_{H_3\cap H_4} = Q_4|_{H_3\cap H_4}$ ($= q$), the pair $(Q_3,Q_4)$ is naturally regarded as a section 
of a line bundle $\mathscr O_{H_3\cup H_4}(2)$.
By the exact sequence 
$$
0 \lra \mathscr O_{\mathbb{CP}^4} \stackrel{\otimes z_3z_4}{\lra} \mathscr O_{\mathbb{CP}^4}(2) \lra \mathscr O_{H_3\cup H_4}(2) \lra 0,
$$
we obtain that the section $(Q_3,Q_4)$ can be extended to a quadratic polynomial  on $\mathbb{CP}^4$, and that such polynomial is unique up to adding a constant multiple of $z_3z_4$.
Explicitly such an extension has to be of the form
\begin{align}\label{Q3}
Q(z_0,z_1,z_2,z_3,z_4) := q(z_0,z_1,z_2) + \sum _{0 \le i\le 3} a_iz_iz_3 + \sum _{0\le i\le 3}b_iz_iz_4
+ c z_3z_4,
\end{align}
where $c$ is an arbitrary constant.
Of course, 
 the hyperquadric satisfies $(Q)\supset \mathscr C_3\cup\mathscr C_4$.
Then if we further impose that  $(Q)$ contains  the point $p_2\in\mathscr C_1\cap\mathscr C_2\subset l=\{z_0=z_1=z_2=0\}$,
then from \eqref{Q3} 
a linear equation for $a_3,b_3, c$ is obtained, from
which $c$ is uniquely determined.
$c$ is real since $a_3$ and $b_3$ are real.
Summarizing up,  we have obtained that once we fix a real quadratic polynomial $q(z_0,z_1,z_2)$, then there exists a unique real quadratic polynomial $Q(z_0,z_1,z_2,z_3,z_4)$ whose restriction to $H_3\cap H_4$ equals $q$
and which goes through the 9 points (1 to 8 and $p_2$).
In particular, if $q=z_0^2-z_1z_2$, then $Q=z_0^2-z_1z_2$.

Next we show that  this polynomial $Q$ (which is uniquely determined from $q$) always contains the double curves $\mathscr C_i$ for any $1\le i\le 5$.
It remains to show $\mathscr C_i\subset (Q)$ for $i=1,2,5$.
For $i=1,2$ this is immediate since $(Q)$ already goes through 5 points on the conic $\mathscr C_i$, which means that it goes through the remaining 1 point $\ol p_2$.
$\mathscr C_5\subset(Q)$ is also immediate if we notice that
as $\mathscr C_i\subset(Q)$ for $1\le i\le 4$,
$(Q)$ already goes through all the 12 points
$\mathscr C_i\cap \mathscr C_5$ for $1\le i\le 4$ on $H_5$ among the 26 points obtained in Section \ref{ss:quadric} and that 
since $h^0(\mathscr O_{H_5\cap Y}(2))=9$, eight points already and
uniquely determine the quadric.
Thus we have proved the existence of a quadratic polynomial $Q$ satisfying $\mathscr C_i\subset (Q)$ for any $1\le i\le 5$.

For the uniqueness  in the sense of the proposition,
let $Q$ and $Q'$ be as stated in the proposition.
Then since both $Q|_{H_3\cap H_4}$ and $Q'|_{H_3\cap H_4}$ are real and
belong to the pencil (determined from the 4 points $\mathscr C_3\cap\mathscr C_4$), 
there exists $(c,c')\in\mathbb R^2$ with $(c,c')\neq(0,0)$
such that $cQ-c'Q'|_{H_3\cap H_4}\in (z_0^2-z_1z_2)$.
Further, the hyperquadric $(cQ-c'Q')$ goes through the points $1$ to $8$ and $p_2$ at least, and  therefore must belong to the ideal $
(z_0^2-z_1z_2)$ by the uniqueness which was already proved.
Thus we get the required uniqueness.
\proofend

\subsection{Defining equation of the branch divisor}\label{ss:defeq}
With the results in the previous 2 subsections, we are ready to provide the main result in this paper:
\begin{theorem}\label{thm:B}
Let $Z$ be any twistor space on $4\mathbb{CP}^2$ containing the surface $S$ (constructed in Section 2) as a real member of $|F|$. 
Let $\Phi_4:Z_4\to Y$ be the generically 2 to 1 morphism canonically obtained from the explicit birational transformations in Section 3, and $B$ the branch divisor of $\Phi_4$.
Let $Q(z_0,z_1,z_2,z_3,z_4)$ be a defining equation of the hyperquadric containing all the 5 double curves, obtained in Proposition \ref{prop:Q}.
Then $B$ is an intersection of the scroll $Y=\{z_0^2=z_1z_2\}$ with the quartic hypersurface defined by
the equation of the form
\begin{align}\label{scrB}
z_0 z_3 z_4 f(z_0,z_1,z_2,z_3,z_4) = 
Q(z_0,z_1,z_2,z_3,z_4)^2 
\end{align}
where $f(z_0,z_1,z_2,z_3,z_4)$ is a linear polynomial with real coefficients.
\end{theorem}

\proof
By Proposition \ref{branchdivisor1} there exists a real hyperquartic such that 
the intersection with  $Y$ is the branch divisor $B$.
Let 
$\mathscr B\subset\mathbb{CP}^4$ be any one of such hyperquartics and
 $F=F(z_0,z_1,z_2,z_3,z_4)$ a defining equation of $\mathscr B$.
We note that $\mathscr B$ is not unique in the sense that $F$ is determined only up to quartic polynomials in the ideal $(z_0^2-z_1z_2)$.
Then for $i=1,2$, the restriction of $\mathscr B$ to a plane $H_i\cap Y=\{z_0=z_i=0\}$ is the twice of the double conic $\mathscr C_i$ (see \eqref{dc3}).
Also, by the choice of $Q$, the restriction of the hyperquadric $(Q)$ to the same plane is $\mathscr C_i$.
These two mean that there exists a constant $c_i\in\mathbb C$ such that $F- c_i Q^2$ belongs to the ideal $(z_0,z_i)$. Namely there exist cubic polynomials $f_i$ and $g_i$ (in $z_0,z_1,z_2,z_3,z_4)$ satisfying
\begin{align}\label{mtp1}
F - c_i Q^2 = z_0 f_i + z_i g_i \quad (i=1,2).
\end{align}
Taking the difference, we obtain 
\begin{align}\label{mtp1'}
(c_1-c_2)Q^2 = z_0 (f_1-f_2) + z_1g_1 - z_2g_2.
\end{align}
If $c_1\neq c_2$, substituting $z_0=z_1=0$, the hyperquadric $(Q)$ restricted to the plane $\{z_0=z_1=0\}$ is defined by $z_2 g_2=0$.
This means that $\mathscr C_1$ is reducible, which cannot happen since it is the image of the twistor line $L_1$
(see \eqref{dc}).
Hence we obtain $c_1=c_2$.
Similarly, for $i=3,4$, considering the restrictions of $F$ and $Q^2$ to the cone $H_i\cap Y=\{z_i=z_0^2-z_1z_2=0\}$, again by coincidence, there exist  a constant $c_i\in\mathbb C$, a cubic polynomial $f_i$ and a quadratic polynomial $g_i$ satisfying
\begin{align}\label{mtp2}
F - c_i Q^2 = z_i f_i + ( z_0^2 - z_1 z_2) g_i \quad (i=3,4).
\end{align}
By \eqref{mtp1} with $i=1$ and \eqref{mtp2}  we obtain
\begin{align}\label{mtp2'}
(c_3-c_1) Q^2 = z_0 f_1 + z_1 g_1- z_i f_i - ( z_0^2 - z_1 z_2)g_i \quad (i=3,4).
\end{align}
From this  we again obtain that the hyperquadric $\{Q=0\}$ restricted to the plane $\{z_0=z_1=0\}$ is 
given by $\{z_i f_i=0\}$, contradicting the irreducibility of $\mathscr C_1$.
Hence we obtain $c_1=c_i$ for $i=3,4$.
Thus we get $c_1=c_2=c_3=c_4$.
If $c_1=0$, by \eqref{mtp1}, we have $F= z_0 f_1 + z_1 g_1$.
But this cannot happen since this means $\mathscr B\supset
\{z_0=z_1=0\}$, contradicting $\mathscr B\cap \{z_0=z_1=0\}=\mathscr C_1$. Hence $c_1\neq 0$. 
So replacing $Q$ with $Q/\sqrt{c_i}$, we may assume that all the four $c_i$-s in \eqref{mtp1} and \eqref{mtp2} are one. 

Next in the expression \eqref{mtp1} we take $f_i$ and $g_i$ in such a way that $g_i$ does not contain $z_0$.
Then since the right hand side of  \eqref{mtp1'} is zero,
$z_1g_1-z_2g_2=0$ follows. 
Hence $g_1\in (z_2)$, and we can write $g_1= z_2h_1$ by a quadratic polynomial $h_1$ which does not contain $z_0$.
Similarly, in the expression \eqref{mtp2} we take $f_i$ and $g_i$ in such a way that $f_3$ and $f_4$ do not belong to the ideal $(z_0^2-z_1z_2)$.
Then this time from \eqref{mtp2} for the case $i=3$ and $i=4$, we obtain
\begin{align}
(z_3f_3-z_4f_4) + (z_0^2 - z_1z_2)(g_3-g_4) = 0.
\end{align}
From the choice of $f_3$ and $f_4$, 
it follows $f_3\in (z_4)$ and $f_4\in (z_3)$.
Hence we can put $f_3=z_4 f_5$ for some quadratic polynomial $f_5$.
From these, we obtain 
\begin{align}\label{mtp4}
F-Q^2 = z_0 f_1 + z_1 z_2 h_1 = z_3 z_4 f_5 + (z_0^2 - z_1 z_2)g_3.
\end{align}
Then since $h_1$ does not contain $z_0$,
from the latter equality we can readily deduce that if we write $f_5 = z_0 f_6 + f_7$
in a way that $f_7$ does not contain $z_0$, then $f_7$ is a multiple of $z_1z_2$, so that 
$f_5 = z_0 f_6 + c z_1 z_2$ for some $c\in\mathbb C$.
Hence by \eqref{mtp4} we obtain 
\begin{align}\label{mtp5}
F- Q^2 = z_3 z_4 (z_0 f_6 + c z_1 z_2)+ (z_0^2 - z_1 z_2)g_3.
\end{align}
Defining a linear polynomial $f_8$ by  $f_6= - cz_0 + f_8$ (so that $f_8$ may contain $z_0$) and substituting into \eqref{mtp5}, we  finally get
\begin{align}\label{mtp6}
F- Q^2 
= z_0z_3z_4f_8 + (z_0^2 - z_1 z_2) ( g_3 - cz_3z_4).
\end{align}
Thus we obtain $F= Q^2 + z_0z_3z_4f_8 + (z_0^2 - z_1 z_2) ( g_3 - cz_3z_4)$.
Hence modulo quartic polynomials in the ideal $(z_0^2-z_1z_2)$, $\mathscr B$ is defined by the equation of the form \eqref{scrB}.
This completes a proof of the theorem. 
\proofend

\vspace{3mm}
From the quartic equation \eqref{scrB} it is immediate to see that the intersection of the hyperquartic $(Q)$ and 
the hyperplanes $H_3$ and $H_4$ are double quadric surfaces, and this is of course consistent with the fact that the restrictions $B|_{H_3\cap Y}$ and $B|_{H_4\cap Y}$ are double curves.
On the other hand, for $i=1,2$, in order to see that $B|_{ H_i\cap Y}$ are also double curves from the equation, we just need to notice that, on the scroll $Y$, $z_i=0$ means $z_0=0$.

In comparison with the case of $3\mathbb{CP}^2$, 
 appearance of the linear polynomial $f$ in our defining 
equation \eqref{scrB} might look strange at first sight.
As the following proposition shows, $f$ comes from the {\em 
fifth} double curve $\mathscr C_5$, which does not exist in
the case of $3\mathbb{CP}^2$:

\begin{proposition}
Up to non-zero constants, the linear polynomial $f$ in \eqref{scrB} is exactly $z_5$ we have defined in \eqref{z5}.
In other words, for the third double quartic curve, we have $\mathscr C_5=\{f=Q=z_0^2-z_1z_2=0\}$.
\end{proposition}

\proof
We will find all hyperplanes in $\mathbb{CP}^4$ (which is the target space of the anticanonical map $\Phi$), which correspond to reducible members of the anticanonical system $|2F|$.
First as above for any hyperplane $H$ defined by the equation of the form $a_0z_0+a_1z_1+a_2z_2=0$, the corresponding 
member $\Phi^{-1}(H)\in |2F|$ is clearly reducible.
(We are including  a non-reduced case.)
Also, $\Phi^{-1}(H_3)$ and $\Phi^{-1}(H_4)$ are reducible since $B\cap H_3$ an $B\cap H_4$ are double curves.
By the same reason, if $H_f=\{f=0\}$, the divisor 
$\Phi^{-1}(H_f)$ is reducible.
Then recalling that  the double quartic curve
 $\mathscr C_5$ is obtained as an image of the third reducible member of $|2F|$ obtained in Proposition \ref{prop:redivisor},
in order to prove the claim of the proposition,
 it suffices to show that 
there exists no other hyperplane $H$ such that 
$\Phi^{-1}(H)$ is  reducible.

If $H$ is such a hyperplane, then either $Y\cap H$ is reducible, or $Y\cap H$ is irreducible (i.e.\,a cone) and $B|_{Y\cap H}$ is a double curve.
The former occurs exactly when $H$ is defined by the equation of the form $a_0z_0+a_1z_1+a_2z_2=0$.
So suppose the latter happens.
Then recalling $B=\mathscr B\cap Y$ and $H|_Y$ is reduced by the assumption, $B|_{Y\cap H}$ can be a double curve only when 
$\mathscr B|_H$ is a double surface.
We  show by algebraic mean that this happens only when $H$ is defined by
one of the 4  factors of the left-hand side of \eqref{scrB}.

If $\mathscr B|_ H$ is a double surface, there exists a quadratic polynomial $q$ on $H$ such that 
$(z_0z_3z_4f-Q^2)|_H=q^2$; namely
\begin{align}\label{aaa}
z_0z_3z_4f\,|_H=(Q|_H)^2 + q^2.
\end{align}
If $H$ is defined by the equation of the form $z_1=b_0z_0+b_2z_2+b_3z_3+b_4z_4$, 
then even after substitution the left-hand side of \eqref{aaa} does not have 
monomial of the form $z_i^3z_j$ for any $1\le i,j\le 4$.
Therefore the coefficient of $z_i^3z_j$ of the right-hand side of \eqref{aaa} must be zero for any $1\le i,j\le 4$.
By an elementary argument, it is possible to show this can happen only when $(Q|_H)^2 + q^2=0$.
Hence by \eqref{aaa} $H$ is equal to one of $H_0, H_3, H_4$ 
and $H_f$.
By symmetry of the equation, if $H$ is defined by a equation of the form 
$z_2=b_0z_0+b_1z_2+b_3z_3+b_4z_4$, 
then  \eqref{aaa} is possible only when $H$ is  
one of $H_0, H_3, H_4$  and $H_f$.
If $H$ is of the form $z_0=b_3z_3+b_4z_4$, 
then the left-hand side of \eqref{aaa} cannot contain
monomials of the form $z_i^4$,
$z_1^3z_i$, $z_2^3z_i$  (for any $i$),
and $z_1^2z_3z_4$, $z_1z_2z_3z_4$, $z_2^2z_3z_4$.
Therefore the coefficients of these monomials of 
the right-hand side of \eqref{aaa} must vanish.
From these, again by an elementary argument it is possible
to show that $(Q|_H)^2 + q^2=0$.
Hence again $H$ has to be one of $H_0, H_3, H_4$  and $H_f$.
The remaining 2 cases immediately follow from symmetry of the equation.
Thus we have shown that 
$\mathscr B|_H$ is a double surface only when 
$H$ is one of $H_0, H_3, H_4$  and $H_f$.
\proofend

\vspace{3mm}
We again emphasize that the role of the 3 double quartic curves is symmetric, and any choice of two leads to  the equation of the form \eqref{scrB}.

\begin{remark}{\em
One may wonder whether the linear polynomial $f$ can be taken as one of the homogenous coordinates on $\mathbb{CP}^4$. 
At least in generic situation this is possible, but 
if we do so, we lose  simplicity of the defining  equation of the scroll $Y$, and it makes  more difficult the counting  the number of effective parameters
in defining equations of $\mathscr B\cap Y$ which will be done  in Section \ref{ss:moduli}.
}
\end{remark}

\subsection{The number  of singularities of the branch divisor }\label{ss:sing}
In this subsection by using the quartic equation obtained in the previous subsection
we  determine the number of singularities 
of the branch divisor of the double covering.
Similarly to the method by Kreussler \cite{Kr89} and
Kreussler-Kurke \cite{KK92},
we resort to {\em topology}\,; more precisely we compute the Euler number of 
the relevant spaces to determine the number of singularities.
Though we require much more complicated computation than the case of $3\mathbb{CP}^2$, we do it since
 this result is crucial for
determining the dimension of the moduli space of the present
twistor spaces.

As in the proof of Proposition \ref{prop:isolated} let $Z_4\stackrel{\mu_5}{\to} Z_5\stackrel{\Phi_5}{\to} Y$ be the Stein factorization of the degree 2 morphism $\Phi_4:Z_4\to Y$.
We already know that 
$\mu_5$ contracts finitely many curves,
whose images are singular points of the branch divisor $B$.
If we put $l_5=\mu_5(l_4)$ and $\ol l_5=\mu_5(\ol l_4)$, 
$Z_5$ also has ordinary double points along $l_5\cup\ol l_5$.
All other singularities of $Z_5$ are lying on singularities of the branch divisor $B$.
Among these singularities we already know that there are 26
singularities listed in Section \ref{ss:quadric},
and 24 points among them are 
 ordinary double points.
The 2 points excluded here are exactly the points 
$\mathscr C_1\cap\mathscr C_2$, for which we still denote by $p_2$ 
and $\ol p_2$ (see Figure \ref{fig:dc} again).
We begin with determining the type of singularities of these 2 points:

\begin{proposition}\label{prop:A3}
At the 2 points $p_2$ and $\ol p_2$, the branch divisor $B$ has $A_3$-singularities.
\end{proposition}

\proof
Recall that $l=\{z_0=z_1=z_2=0\}$, and $\mathscr C_1\cap\mathscr C_2=B\cap l=\{Q=0\}\cap l$.
As above let $p_2$ be any one of the 2 points and we work in a neighborhood of $p_2$.
We put $x:=z_0/z_4, y:= z_1/z_4, z:= z_2/z_4$ and $u:=Q$.
Then by transversality for the intersection  of $Q$ and $l$,
we can use $(x,y,z,u)$ as coordinates in a neighborhood of $p_2$ in $ \mathbb{CP}^4$,
and noticing $z_3z_4f\neq 0$ at $p_2$,
we may suppose that the hyperquartic \eqref{scrB} is defined by a very simple equation, $x=u^2$.
Since $Y$ is defined by $x^2=yz$, we deduce that 
 $p_2$ is an  A$_3$-singular point of the surface $B$.
By  reality, $\ol p_2$ is also an A$_3$-singular point.
\proofend

\vsp
Next, as the transformation from $Z$ to $Z_4$ is explicit, it is easy to show the following:

\begin{proposition}\label{euler4}
For the  variety $Z_4$ we have $e(Z_4)=10$.
\end{proposition}

\proof
Since $Z$ is a twistor space on $4\mathbb{CP}^2$, we have $e(Z)= 2+2(b_2(4\mathbb{CP}^2)+1)= 12$.
Because  the blowup $\mu_1$ replaces two disjoint $\mathbb{CP}^1$-s by two $\mathbb{CP}^1\times\mathbb{CP}^1$-s, we have $e(Z_1)=12 + 4 = 16$.
Then since a flop does not change the Euler number 
we obtain $e(Z_3)= 16$.
Finally looking Figure \ref{fig:trans} (d),  the exceptional divisor
$E_1\cup\ol E_1$ of the contraction $\mu_4:Z_3\to Z_4$ has Euler number 8, and 
the image $l_4\cup \ol l_4$ of the exceptional divisor has Euler number 2.
Hence we obtain $e(Z_4)= 16 - (8-2) = 10$.
\proofend

\vsp
The next result means that the 26 points that we have already found are {\em not} all singularities of the branch divisor $B$, but  in the generic situation $B$ has extra 6 ordinary double points:

\begin{theorem}\label{thm:6sing}
Let $\{b_1,\cdots,b_k\}$ be the set of all singular points of $B$ which are different from
the 26 singular points listed in Section \ref{ss:quadric}.
Let $\mu_i$ be the Milnor number of the singular point $b_i$, and put
$\beta_i:=\mu_5^{-1}(b_i)$ for the exceptional curve of $\mu_5$ over the point $b_i$.
(Of course we do not assume irreducibility of $\beta_i$.)
Then we have the relation
\begin{align}\label{strange}
\sum_{i=1}^k \left\{e(\beta_i)+\mu_i -1\right\} = 12.
\end{align}
In particular, if all the singularities are ordinary double points, we have $k=6$.
\end{theorem}

\proof
First since $\mu_5:Z_4 \to Z_5$ replaces each of the 24 ordinary double points by smooth $\mathbb{CP}^1$ and also replaces the singular point $b_i$ by the curve $\beta_i$ for $1\le i\le k$, we obtain
\begin{align}
e(Z_4) = e(Z_5) + 24 + \sum_{1\le i\le k} \{ e(\beta_i)-1 \}.
\end{align}
On the other hand by the double covering $Z_5\to Y$ we have
$
e(Z_5) = 2 e(Y) - e(B).
$
Further as $Y$ is obtained from the $\mathbb{CP}^2$-bundle $\tilde Y$ over $\mathbb{CP}^1$, and it replaces $\mathbb{CP}^1\times\mathbb{CP}^1$ by $l\simeq\mathbb{CP}^1$, we have $e( Y) = e(\tilde Y) - (4-2) = 6 -2= 4$.
Hence we have $e(Z_5) = 8 - e(B)$, giving
\begin{align}\label{euler7}
e(Z_4) = 32 - e (B) + \sum_{1\le i\le k} \{ e(\beta_i ) -1 \}.
\end{align}
Next for computing $e(B)$ let $D$ be a general member of the system $|\mathscr O_Y(4)|$.
Then since the scroll $Y$ has ordinary double points along the line $l$, the divisor $D$ has ordinary double points at the 4 points $D\cap l$.
As before let $\nu:\tilde Y\to Y$ be the blowup at $l$, and let $\tilde D$ be the strict transform of $D$.
By Bertini's theorem we may suppose that $\tilde D$ is non-singular.
We shall compute $e(\tilde D)$.

As before write $\mathfrak f:=\tilde{\pi}^*\mathscr O_{\Lambda}(1)\in H^2(\tilde Y,\mathbb Z)$.
From the standard relationship of the total Chern class $c(T_{\tilde Y}|_{\tilde D})= c(T_{\tilde D})\cdot
c(N_{\tilde D/\tilde Y})$ 
and adjunction formula, we readily obtain 
\begin{align}\label{euler9}
e(\tilde D) = c_2(T_{\tilde D}) = c_2(T_{\tilde Y})\cdot \tilde D + (K_{\tilde Y}+ \tilde D)\cdot \tilde D\cdot \tilde D,
\end{align}
where the dot means the product in $H^*(\tilde Y,\mathbb Z)$.
As in the proof of Proposition \ref{branchdivisor1} let $\mathscr O(0,1):= (\nu^*\mathscr O_l(1))|_{\Sigma}$, so that 
$\mathscr O(1,0)= \mathfrak f|_{\Sigma}$.
Then by using $N_{\Sigma/\tilde Y}\simeq \mathscr O(-2,1)$ and the adjunction formula applied to 
a fiber of $\tilde{\pi}$, we readily obtain $K_{\tilde Y}\sim -3\Sigma - 6 \mathfrak f$.
Further,  in the cohomology ring of $\tilde Y$, we have $\Sigma^3 = N_{\Sigma/\tilde Y}
\cdot N_{\Sigma/\tilde Y} = -4$,
$\Sigma^2\cdot\mathfrak f = N_{\Sigma/\tilde Y}\cdot \mathfrak f = \mathscr O(-2,1) \cdot \mathscr O(1,0) = 1$,
 $\Sigma\cdot \mathfrak f^2 = 0$ and $\mathfrak f^3 = 0$.
 Furthermore by recalling $\nu^*\mathscr O(1)\sim \Sigma +2\mathfrak f$
(see \eqref{hs22}), we obtain
 $K_{\tilde Y} + \tilde D\sim (-3\Sigma -6\mathfrak f) + 4(\Sigma +2\mathfrak f)=\Sigma +2\mathfrak f$.
From these we readily get $(K_{\tilde Y}+ \tilde D)\cdot \tilde D\cdot \tilde D = 32$.

For computing $c_2(T_{\tilde Y})\in H^4(\tilde Y,\mathbb Z)$, as  generators of $H^4(\tilde Y,\mathbb Z)$ we take any element $\zeta\in |\mathscr O_{\Sigma}(1,0)|$ and $\eta\in |\mathscr O_{\Sigma}(0,1)|$, viewed
as submanifolds in $\tilde Y$, and 
put $c_2(T_{\tilde Y}) = a \zeta + b \eta$.
From the exact sequence associated to the inclusion $\Sigma\subset \tilde Y$,
 we immediately obtain $c_2(T_{\tilde Y})|_{\Sigma} = c_1(\Sigma) \cdot c_1(N_{\Sigma/\tilde Y})
 + c_2(\Sigma)$.
 Then since $c_1(\Sigma) = \mathscr O(2,2),\, c_1(N_{\Sigma/\tilde Y})= (-2,1)$ and 
 $c_2(\Sigma) = e(\Sigma) = 4$, 
 we obtain $c_2(T_{\tilde Y})|_{\Sigma} = 2$.
 On the other hand, from the inclusion $\mathfrak f\subset\tilde Y$ we readily obtain 
 $c_2(T_{\tilde Y}) |_{\mathfrak f} = 3$.
 Further, in the cohomology ring of $\tilde Y$ we have $\zeta\cdot\Sigma =\mathscr O(1,0)\cdot
 \mathscr O(-2,1) = 1, \, \eta\cdot\Sigma = \mathscr O(0,1) \cdot \mathscr O(-2,1) = -2,
 \, \zeta\cdot \mathfrak f = \mathscr O(1,0)\cdot\mathscr O(1,0) = 0$
 and $\eta\cdot \mathfrak f = \mathscr O(0,1)\cdot\mathscr O(1,0) = 1$.
 Therefore by restricting to $\Sigma$ and $\mathfrak f$ respectively, we get
 $a-2b =2 $ and $b=3$. 
 Hence $a=8$.
 Therefore we obtain 
$
 c_2(\tilde Y)\cdot \tilde D = (8 \zeta + 3 \eta) \cdot 4 ( \Sigma + 2 \mathfrak f ),
$
which is readily computed to be $32$.
Hence from \eqref{euler9} we obtain
$ e (\tilde D) = 32 + 32 = 64$.

As $\tilde D\to D$ contracts four $\mathbb{CP}^1$-s to 4 points, we have $e(D)=60$.
Then $D$ is obtained from  the actual branch divisor $B$ by 
(a) smoothing the $24$ nodes, (b) smoothing $k$ singular points $b_1,\cdots,b_k$,
and (c) deforming each of the two $A_3$-singularities (which is exactly $\mathscr C_1\cap \mathscr C_2=\{p_2,\ol p_2\}$) 
to two $A_1$ singularities.
Adding the Milnor number of the singularities for the cases (a) and (b), 
and also taking the difference of the Milnor number of $A_3$-singularity and two $A_1$-singularities into account,
we obtain
\begin{align}
e ( B ) = e ( D ) - 26 - \sum_{1\le i\le k} \mu_i,
\end{align}
and hence $e(B) = 34 - \sum_{1\le i\le k} \mu_i$.
Substituting this into \eqref{euler7} and using Proposition \ref{euler4}, we obtain \eqref{strange}.
\proofend

\begin{remark}{\em
If the blown-up 8 points on $\mathbb{CP}^1\times\mathbb{CP}^1$ are arranged as in Figure \ref{fig:square} is in a general position 
(in certain precise sense), then $C_2$ and $\ol C_2$ are all curves that are contracted to points by the bi-anticanonical map.
But if  the 8 points are in a special position (in certain precise sense), then the map  contracts extra curves.
It is not difficult to classify all positions which yield
this situation.    
The appearance of this kind of curves is exactly the reason why the anticanonical map of the 
twistor spaces contracts some  rational curves which cannot be found from the equation of the branch divisor.
}
\end{remark}

\section{Moduli space and genericity of the twistor spaces}
\label{s:moduli}

\subsection{Dimension of the moduli space}
\label{ss:moduli}
In this subsection we compute the dimension of the moduli space of our twistor spaces by counting the number of effective parameters, and verify that it agrees with the dimension of the cohomology group which is relevant to the 
present case.

For the former purpose, 
we recall from Section \ref{s:anticanonical} that $Z$ canonically determines a birational model $Z_4$ and the degree 2 morphism $\Phi_4:Z_4\to Y$, 
and from Section \ref{ss:defeq} the branch divisor of $\Phi_4$ is an intersection of $Y$ with the quartic surface defined by
\begin{align}\label{scrB2}
z_0 z_3 z_4 f(z_0,z_1,z_2,z_3,z_4) = 
Q(z_0,z_1,z_2,z_3,z_4)^2, 
\end{align}
where $f$ and $Q$ are linear and quadratic polynomials with real 
coefficients respectively.
In this subsection we denote this quartic hypersurface by $\mathscr B(f,Q)$.
Since the quartic hypersurface $\mathscr B(f,Q)$ uniquely determines the double cover via the natural quadratic map
$\mathscr O_Y(2)\to \mathscr O_Y(4)$ (which takes squares), up to small resolutions,  $Z$ is uniquely determined by the quartic hypersurface.
Of course, $f$ has 5 coefficients and $Q$ has 15 coefficients,
so the equation \eqref{scrB2} contains $20$ parameters.
Further it is elementary to see that  two pairs $(f,Q)$ and 
$(f',Q')$ of linear and quadratic polynomials (over $\mathbb R$) determine 
the same hyperquartic surface if and only if $(f',Q')=(c^2f,c\,Q)$ for some $c\in\mathbb R$.
This decreases the number of parameters by one.
On the other hand, projective transformations which preserve $Y$ and which preserves {\em the form of} the equation \eqref{scrB2} have to be $(z_0,z_1,z_2,z_3,z_4)\mapsto (abz_0,az_1,bz_2,cz_3,dz_4)$ for some $a,b,c,d\in\mathbb R^{\times}$.
(Here we are only considering transformations which are
homotopic to the identity.)
If a pair $(f',Q')$ is obtained from a pair
$(f,Q)$ by one of these transformations, then the intersections 
$Y\cap \mathscr B(f,Q)$ and $Y\cap \mathscr B(f',Q') $
are mutually biholomorphic, so that they define mutually isomorphic double cover.
But taking an effect of the above equivalence $(c^2 f,c\,Q) = (f, Q)$ into account,
we can suppose $ab=1$ and therefore these projective transformations
decrease the number of parameters by 3.
Thus up to now the number of parameters is $20-(1+3)=16$.
However what we have to consider is not the hyperquartics \eqref{scrB2} themselves 
but the intersection with $Y$;
namely if $Q'=Q+c(z_0^2-z_1z_2)$  for some
$c\in\mathbb R$, then we have the coincidence
$Y\cap \mathscr B(f,Q)= Y\cap \mathscr B(f,Q')\subset\mathbb{CP}^4$.
Clearly these transformations are not included in the above projective
transformations,
so they drop the dimension by one.
Thus we have obtained that the number of effective parameters
in the quadratic hypersurface \eqref{scrB2} is 15.
Finally, by Theorem \ref{thm:6sing},  the pair $(f,Q)$ must satisfy the constraint that $\mathscr B(f,Q)\cap Y$ has extra 6 ordinary double points in general, which decreases the number of  parameters by 6.
Therefore we conclude that {\em the space of  isomorphic classes of the divisors of the form $\mathscr B(f,Q)\cap Y$, which can be the branch divisor for the twistor spaces under consideration, is $15$-dimensional}.

Next we compute the dimension of the moduli space of our twistor spaces by determining the dimension of the first cohomology group of
an appropriate subsheaf of the tangent sheaf.
We begin with a computation for the full moduli space.

\begin{proposition}\label{prop:13dim}
Let $Z$ be a twistor space on $4\mathbb{CP}^2$ which contains the surface $S$ constructed in Section \ref{s:1} as a real member of $|F|$. Then we have $H^i(Z,\Theta_Z)=0$ for $i\neq 1$ and 
$h^1(Z,\Theta_Z) = 13$.
\end{proposition}

\proof
As computed in \cite{LB92}, for any twistor space on 
$n\mathbb{CP}^2$,
by the Riemann-Roch formula, we have $\chi(\Theta_Z) = 15 - 7n$.
Also, since $Z$ is Moishezon and $|F|$ has an irreducible member, we have $H^2(\Theta_Z)=0$ by \cite{C91}.
Further we always have $H^3(\Theta_Z)=0$. 
Hence it suffices to show $H^0(\Theta_Z)=0$.
Let $\Aut_0Z$ be the identity component of the holomorphic automorphism group of $Z$.
Then the real part $(\Aut_0Z)^{\sigma}$ is naturally identified with the identity component of
conformal automorphism group of the self-dual structure. 
Also, since $|F|$ has just 2 irreducible components,
$(\Aut_0Z)^{\sigma}$ acts on $S_1^+\cup S_1^-$.
Hence as the twistor projections $S_1^+\to 4\mathbb{CP}^2$ and $S_1^-\to 4\mathbb{CP}^2$ are 
of degree 1, $(\Aut_0Z)^{\sigma}$ acts effectively
on $S_1^+\cup S_1^-$.
Furthermore, the degree 1 divisor $S_1^+$ 
and $S_1^-$ are obtained from 
$\mathbb{CP}^2$ by blowing-up 4 points, exactly 3 of which are collinear. From this it readily follows that 
the subgroup of $\Aut\, S_0^+$ which consists of automorphism preserving the twistor line
$L_1$ is 0-dimensional. Hence so is $\Aut (S_0^+\cup S_0^-)$.
Thus $\Aut_0Z$ cannot be of positive dimension.
Therefore $H^0(\Theta_Z)=0$.
\proofend

\vsp
From the proposition, the real part of the Kuranishi family of our twistor space $Z$ is 13-dimensional.
Of course, generic twistor spaces on $4\mathbb{CP}^2$ is algebraic dimension 1 and   generic members of the Kuranishi family have the same property.
In order to restrict to the Moishezon twistor spaces under consideration, we show the following.

\begin{proposition}\label{prop:dfmpair}
Let $Z$ and $S$ be as in Proposition \ref{prop:13dim} and 
$C_1$ and $\ol C_1$ the $(-3)$-curves on $S$. 
Then deformation theory of the pair $(Z,C_1\sqcup \ol C_1)$ is unobstructed and its  Kuranishi family is 9-dimensional.
Further, for all sufficiently small deformations preserving the
real structure in the Kuranishi family, 
the twistor spaces contain a non-singular surface constructed in Section \ref{s:1}  as a real member of $|F|$. 
\end{proposition}

Of course, the last property means that the deformed spaces
are still the Moishezon twistor spaces under consideration.

\noindent
{\em Proof of Proposition \ref{prop:dfmpair}.}
Let $\Theta_{Z,\, C_1+\ol C_1}$ be the subsheaf of $\Theta_Z$ whose germs are vector fields that are tangents to $C_1$ and $\ol C_1$.
For the former claim on the Kuranishi family, it suffices to show that 
$H^2(Z,\Theta_{Z,\, C_1+\ol C_1})=0$
and $h^1(Z,\Theta_{Z,\, C_1+\ol C_1})=9$.
Recalling $N_{C_1/Z}\simeq N_{\ol C_1/Z}\simeq\mathscr O(-2)^{\oplus 2}$, we obtain the standard exact sequence
$$
0 \lra \Theta_{Z, C_1+ \ol C_1} \lra \Theta_Z \lra \mathscr O_{C_1}(-2)^{\oplus 2} \oplus \mathscr O_{\ol C_1}(-2)^{\oplus 2} \lra 0,
$$
which induces an exact sequence
\begin{align}
0 \lra H^1(\Theta_{Z,C_1+\ol C_1}) \lra
H^1(\Theta_Z) \lra \mathbb C^4 \lra
 H^2(\Theta_{Z,C_1+\ol C_1}) \lra 0.
\end{align}
Hence with the aid of Proposition \ref{prop:13dim} we have only to show $H^2(\Theta_{Z,C_1+\ol C_1})=0$.
For this we first deduce from duality and rationality that
$H^2(\Theta_S(-C_1-\ol C_1))=0$, which implies, 
from the exact sequence
$0 \to \Theta_S(-C_1-\ol C_1) \to \Theta_{S,C_1+\ol C_1}
\to \Theta_{C_1\sqcup\ol C_1} \to 0 $, that
 $H^2(\Theta_{S,
C_1+\ol C_1})=0.$
Moreover, noting $N_{S/Z}\simeq F|_S\simeq -K_S$, 
we have an exact sequence 
$$
0 \lra \Theta_{S, C_1+\ol C_1} \lra 
\Theta_{Z, C_1+\ol C_1}|_S \lra 
-K_S\otimes \mathscr O_S(-C_1-\ol C_1) \lra 0.
$$
For the last term we have $-K_S\otimes \mathscr O_S(-C_1-\ol C_1) \simeq \mathscr O_S (C_2 +\ol C_2) $, and it is easy to see
$H^2(\mathscr O_S(C_2 + \ol C_2)) = 0$.
Hence for the middle term we obtain $H^2 (\Theta_{Z, C_1+\ol C_1}|_S) =0$.
Then by the exact sequence
$
 0 \lra \Theta_Z(-S) \lra 
\Theta_{Z, C_1+\ol C_1} \lra 
\Theta_{Z, C_1+\ol C_1}|_S \lra 0
$ and $H^2(\Theta_Z(-S))=0$ \cite{C91}, we finally obtain $H^2(\Theta_{Z,C_1+\ol C_1}) = 0$,
as claimed.

For the latter claim about the existence of the surface in the deformed space, let $Z_t$ be any one of the deformed twistor space which is sufficiently close to the original $Z$, and $C_{1t},\ol C_{1t}\subset Z_t$ the curves corresponding to the original curves $C_1$ and $\ol C_1$.
Let $F_t$ be the fundamental line bundle on $Z_t$.
Then as $\dim |F|=1$, we may suppose $\dim |F_t| = 1$ 
 by upper-semicontinuity of dimensions of cohomology groups under  deformations and the Riemann-Roch formula
$ \chi(F_t) = 2.$
We also have an invariance $F_t\cdot C_{1t} = F\cdot C_1$, and the latter is equal to $K_S^{-1}\cdot C_1 =-1$, and therefore
$F_t\cdot C_{1t} = -1$.
This means that $C_{1t}$ and $\ol C_{1t}$ are base curves of the pencil $|F_t|$.
Let $S_t\in |F|$ be any real irreducible member.
Through the Kuranishi family, this surface can be regarded as 
a small deformation of some real irreducible $S\in |F|$,
which means that $S_t$ is obtained from $\mathbb{CP}^1\times\mathbb{CP}^1$ by moving the blownup 8 points from the original positions (indicated as in Figure \ref{fig:square}).
But since $S_t$ contains the curves $C_{1t}$ and $\ol C_{1t}$ as $(-3)$-curves, the property that  3 points belong to a $(1,0)$-curve must be preserved.
This means that the structure of  $S_t$ is the same as that of the original $S$, and we are done.
\proofend

\subsection{Genericity of the twistor spaces}
\label{ss:gen}
In this subsection, by using a theorem of Pedersen-Poon about structure of real irreducible members of $|F|$, we show that the present twistor spaces are in a sense generic among all Moishezon twistor spaces on $4\mathbb{CP}^2$.
We first recall the theorem of Pedersen-Poon \cite{PP94} in a precise form:

\begin{proposition}\label{prop:PP94}
Let $Z$ be a twistor space on $n\mathbb{CP}^2$ and $S\in |F|$  a real irreducible member.
Then $S$ is non-singular with $K_S^2=8-2n$, and the set of twistor lines lying on $S$ is exactly the real part of a real pencil whose self-intersection number is zero.
Moreover there is a birational morphism $\epsilon:S\to\mathbb{CP}^1\times\mathbb{CP}^1$ preserving the real structure, such that  the twistor lines are mapped to $(1,0)$-curves.
\end{proposition}

\begin{figure}
\includegraphics{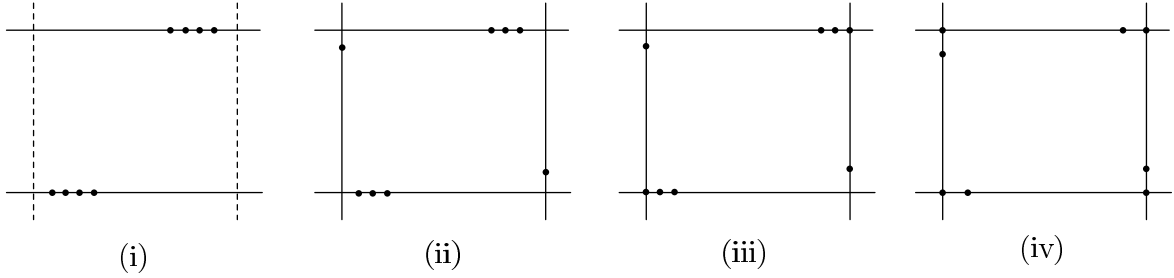}
\caption{Possible configurations of distinct 8 points 
on $\mathbb{CP}^1\times\mathbb{CP}^1$ for Moishezon twistor spaces, except for the Campana-Kreussler's case. }
\label{fig:square2}
\end{figure}

Thus $S$ is always obtained from $\mathbb{CP}^1\times\mathbb{CP}^1$ by
blowing up $2n$ points, where some of the points might be infinitely near in general.
As is well-known 
the position of the blowing up points 
has a strong effect on algebraic structure of twistor spaces.
For example, if a twistor space $Z$ contains $S$ that is obtained from  the $2n$ points  lying on an irreducible $(1,2)$-curve, then it follows $\dim |F|=2$, and detailed structure of such twistor spaces is investigated by Campana-Kreussler \cite{CK98}. 
Then in terms of the configuration of the blowing up points, the genericity of our twistor spaces refers  the following property:

\begin{proposition}\label{prop:density}
Let $Z$ be a Moishezon twistor space on $4\mathbb{CP}^2$ which is not of Campana-Kreussler type.
Suppose that there exists a real irreducible member
$S\in |F|$ such that the images of the 8 exceptional curves
of the blowing-down $\epsilon:S\to\mathbb{CP}^1\times\mathbb{CP}^1$ in Proposition \ref{prop:PP94} 
can be taken as distinct points.
Then the configuration of the 8 points falls into exactly one of Figure \ref{fig:square2}.
\end{proposition}

%
For the proof, we first show the following.

\begin{proposition}\label{prop:gen1}
If $Z$ is a twistor space on $4\mathbb{CP}^2$ which satisfies $\dim |F|=2$,
then $Z$ is either a Campana-Kreussler twistor space,
or otherwise non-Moishezon.
\end{proposition}

\proof
Let $S\in |F|$ be a real irreducible  member, which is necessarily 
non-singular as above.
By the assumption, $S$ satisfies $\dim |K_S^{-1}|=1$.
Let $\epsilon:S\to \mathbb{CP}^1\times\mathbb{CP}^1$ be the 
birational morphism fulfilling the properties of Proposition \ref{prop:PP94}.
(The images of the exceptional curves of $\epsilon$ can be infinitely near.)
Since $\epsilon$ is a composition of blowdowns, 
by the canonical bundle formula for blowups, the image of 
the pencil $|K_S^{-1}|$ by $\epsilon$ necessarily has to be  a pencil of anticanonical curves on  $\mathbb{CP}^1\times\mathbb{CP}^1$; namely a pencil of $(2,2)$-curves.
Let  $\mathscr  P$ be this pencil on  $\mathbb{CP}^1\times\mathbb{CP}^1$.
Then again by the canonical bundle formula 
all the images of the exceptional curves of $\epsilon$
must be contained in the base locus of $\mathscr P$.

Suppose that general members of the pencil $\mathscr P$ are irreducible.
Then the pencil $|K_S^{-1}|$ does not have a fixed component.
If this pencil has a base point, by taking a sequence of blowups $\tilde S\to S$, we obtain a morphism $\tilde S\to\mathbb{CP}^1$, which is, again by the canonical bundle formula,  necessarily the anticanonical map
on $\tilde S$.
Therefore, the morphism $\tilde S\to\mathbb{CP}^1$ must be an elliptic fibration.
But then by the canonical bundle formula for elliptic surfaces, 
we obtain $c_1^2(\tilde S) = 0$.
Since $c_1^2(S)=0$, this means that $\tilde S$ and $S$ are biholomorphic. Hence the pencil $|K_{S}^{-1}|$ is base point free, and the anticanonical map induces an elliptic fibration 
$S\to \mathbb{CP}^1$.
This implies the anti-Kodaira dimension of $S$ is one, which means that $Z$ is non-Moishezon.

So in the sequel we suppose that general members of the pencil $\mathscr P$ on $\mathbb{CP}^1\times\mathbb{CP}^1$ are reducible.
Then if $\mathscr P$  does not have
a fixed component, we have  $\dim|K_S^{-1}|\ge 2$, which contradicts our assumption.
Hence $\mathscr P$ has a fixed component.
Let $C_0$ be any one of its irreducible components.
Then   among the image points of the exceptional curves of $\epsilon$,
there exists at least 1 point on $C_0$
because otherwise $C_0$ is not a fixed component.
Suppose that $C_0\in |\mathscr O(0,1)|$.
Then $\ol C_0\neq C_0$ by the induced real structure on
$\mathbb{CP}^1\times\mathbb{CP}^1$, and $\ol C_0$ is also a 
fixed component of $\mathscr P$.
Hence the movable part of $\mathscr P$ must be a free 1-dimensional subsystem of $|\mathscr O(2,0)|$,
or the system $|\mathscr O(1,0)|$ itself with another fixed $(1,0)$-curve $C'_0$.
But the former  cannot occur because general members of the movable part of the 1-dimensional subsystem would be
reducible by freeness and both components actually move, so that 
all the images of the exceptional curves of $\epsilon$ have to be contained in 
$C_0\cup\ol C_0$, which means $\dim |K_S^{-1}| = 2$.
So suppose the latter is the case.
Then the fixed $(1,0)$-curve $C'_0$
cannot be real, since if so, we would have $C'_0=\epsilon(L)$ for some twistor line $L\subset S$ by the property of $\epsilon$, whereas  on $C'_0$ there is at least one point among the images of the exceptional curves of $\epsilon$, which means $L^2<0$ on $S$.
Hence $\ol C'_0\neq C'_0$.
But this is impossible since $C_0+\ol C_0 + C'_0 + \ol C'_0$, which is clearly $(2,2)$-curves, would be  fixed components of 
the pencil $\mathscr P$ of $(2,2)$-curves.
Thus we obtained $C_0\not\in |\mathscr O(0,1)|$;
namely $\mathscr P$ does not have a $(0,1)$-curve as a fixed component.

Also if a fixed component $C_0$ is a $(1,0)$-curve, then it cannot be real by the same reason.
Hence the movable part of $\mathscr P$ is a free 1-dimensional subsystem of $|\mathscr O(0,2)|$, or  the system $|\mathscr O(0,1)|$ with another fixed component $C'_0\in|\mathscr O(0,1)|$.
But the former implies $\dim |K_{S}^{-1}| = 2$ by the same argument as above, and the latter cannot occur since this time there is no real $(0,1)$-curve.
Thus the fixed  component $C_0$ of $\mathscr P$ cannot be a $(1,0)$-curve.
Further we have $C_0\not\in |\mathscr O(1,1)|$, since any $(1,1)$-curve is not real and hence $\ol C_0\neq C_0$ has to be also a base curve, which contradicts that $\mathscr P$ is a pencil.
Similarly we have $C_0\not\in |\mathscr O(2,1)|$ by the real structure.
Hence $C_0$ must be in the remaining possibility, $C_0\in | \mathscr O(1,2)|$.
In this case $\dim |K_S^{-1}| = 1$ means that all the images of 
the exceptional curves of $\epsilon$ belong to $C_0$.
This implies that the structure of $S$ is exactly as  in the case of Campana-Kreussler,
and we are done.
\proofend

\vsp\noindent
{\em Proof of Proposition \ref{prop:density}}.
First by a result by Kreussler \cite[Theorem 6.5]{Kr98},
on $n\mathbb{CP}^2$ with $n\ge 3$ we always have $\dim |F|\le 3$ and the equality holds iff $Z$ is a LeBrun twistor space \cite{LB91}.
Suppose that $Z$ is a LeBrun twistor space.
Then it is well-known that
a configuration of $2n$ points for generic real
irreducible member $S\in |F|$ is as in (i) of 
Figure \ref{fig:square2}.


So suppose that $Z$ is a Moishezon twistor space on $4\mathbb{CP}^2$ which is different from LeBrun's nor Campana-Kreussler's, and let $S$ be a real irreducible member of $|F|$ such that the $2n$ points on $\mathbb{CP}^1\times\mathbb{CP}^1$ are distinct.
By Proposition \ref{prop:gen1} we have $\dim |F|=1$.
This means $\dim | K_S^{-1} | = 0$.
Let $C$ be the unique anticanonical curve,
 $\epsilon : S\to \mathbb{CP}^1\times\mathbb{CP}^1$ the birational morphism as in Proposition \ref{prop:PP94} whose
images of the exceptional curves are distinct. 
We
put $C_0:=\epsilon(C)$, which is necessarily a real $(2,2)$-curve.
All the 8 points are on $C_0$.
If $C_0$ is irreducible,  $C_0$ must be a non-singular elliptic curve by using the real structure, and 
from this we readily see that $h^0(mK_S^{-1})\le m$ for all $m> 0$, which means that $Z$ is non-Moishezon.
Hence $C_0$ is reducible.
Taking the form of  the induced real structure on $ \mathbb{CP}^1\times\mathbb{CP}^1$ into account,
we can easily show that the decomposition of $C_0$ into irreducible components is one of the following 3 types:
(a) $(1,0) + (0,1) + (1,0) + (0,1)$, 
(b) $(1,1) + (1,1)$, both components being irreducible, or
(c) $(1,2) + (1,0)$, both components being irreducible.
We note that since $h^0(K_S^{-1})=1$, on any
of these components, there exists at least 1 point among
the 8 points.
Repeating an argument  in the last part of the proof of Proposition \ref{prop:gen1}, we deduce that (c) cannot happen
under our assumption.
If $C_0$ has a multiple component, 
 since there exists no  real $(0,1)$-curve, 
it must be  a real $(1,0)$-curve.
But this cannot happen since as remarked above among
the 8 points there is at least one point on any irreducible component of $C_0$, contradicting the family of twistor lines
on $S$.
Hence in both cases (a) and (b) $C_0$ has no multiple component.

Next we show that in the case (a) there is an irreducible component of $C_0$ on which precisely $ 3 $ points among the 8 points  lie.
If not,
then because we are excluding LeBrun twistor spaces,
 on each of the 4 irreducible components exactly 2 points are lying among the 8 points.
In this case, the restriction $K_S^{-1}|_{C}\simeq [C]|_{C}$ belongs to Pic\,$^0C\simeq\mathbb C^*$, which again implies
 $h^0(-mK_S)\le m$ for any $m> 0$ as in the above case.
This implies that $Z$ is not Moishezon.
Hence the component actually exists.

Next we prove that if $C_0$ is in the case (b), there exists another birational morphism $\epsilon':S\to 
 \mathbb{CP}^1\times\mathbb{CP}^1$ preserving the real structure such that $C_0'=\epsilon'(C)$ 
 falls into the case (a), and such that the images of twistor lines in $S$ are $(1,0)$-curves.
For this we  write $C_0=C_1+\ol C_1$ with $C_1$ and $\ol C_1$
being irreducible $(1,1)$-curves.
If exactly 4 points belong to $C_1$, then
the remaining 4 points belong to $\ol C_1$,
and also no point coincides with the 2 points $C_1\cap \ol C_1$.
Then by a similar reason for the case (a), this implies that $Z$ is not Moishezon.
So the 2 points $C_1\cap \ol C_1$ are included in the 8 points.
Therefore $\epsilon$ factors as $S\to S_1\to  \mathbb{CP}^1\times\mathbb{CP}^1$
where the latter arrow is the blowup at $C_1\cap \ol C_1$.
We obtain  6 points on $S_1$ as
the images of the exceptional curves of $S\to S_1$.
These 6 points are not on the exceptional curves of $S_1\to
\mathbb{CP}^1\times\mathbb{CP}^1$ by the assumption
that the 8 points are distinct.
This implies that 
 the unique anticanonical curve
$C$ on $S$ is a cycle of 4 rational curves, whose self-intersection  numbers are 
 $(-3), (-1), (-3), (-1)$.
Then for another blowdown $\epsilon'$ in the proposition,
it is enough to choose a blowingdown $S_1\to \mathbb{CP}^1\times\mathbb{CP}^1$ which does not contract the 2 exceptional
curves of the original $S_1\to\mathbb{CP}^1\times\mathbb{CP}^1$, and composite it with the morphism $S\to S_1$.
Thus we obtained another $\epsilon'$ as claimed, and 
we can neglect the case (b).

Hence $C_0$ can be supposed to be in the case (a) and that there is at least one component on which exactly 3 points among 8 points lie.
This directly means the 8 points have to be put on $C_0$ arranged as in (ii), (iii) or (iv), as claimed.  
\proofend

\vsp
Needless to say, the case (ii) of 
Proposition \ref{prop:density} is exactly the 
situation we have investigated in this paper.
Since it is clear that 
all other 3 cases 
((i), (iii) and (iv)) can be  obtained as  small deformations
of the case (ii), it would be reasonable to say that among the surface $S$ obtained from the 8 points arranged as in (i)--(iv), the case (ii) is most generic.
By deformation theory including a co-stability theorem
of Horikawa \cite{Hor76}, the same is true for 
the twistor spaces containing these surfaces.
Namely any twistor spaces on $4\mathbb{CP}^2$ which has $S$
obtained from (i), (iii) and (iv) as a real member of $|F|$
can be obtained as a limit of the twistor spaces investigated
in this paper.
In particular, the present twistor spaces can be obtained
as a small deformation of a LeBrun twistor space,
and this proves the existence of our twistor spaces.
We also remark that by using the Horikawa's theorem,
it is possible to show that the present twistor spaces
can also be obtained as a small deformation of the twistor spaces studied in \cite{Hon-II} (on $4\mathbb{CP}^2$, of course).
These are the reason why we call the present twistor
spaces to be generic.

Finally we remark that a converse of Proposition \ref{prop:density}
 also holds. Namely
if a twistor space $Z$ on $4\mathbb{CP}^2$ has real irreducible $S\in |F|$ which is obtained from the 8 points in the case (ii), (iii), or (iv),
then  $\dim |F|=1$ and $Z$ is Moishezon. 
Concerning  structure of these twistor spaces, the case (iii) 
can be regarded as a mild degeneration of the present 
twistor spaces, in the sense that 
the twistor space still has a double covering structure over the scroll $Y$ by the anticanonical system.
These twistor spaces (having $S$ obtained from the configuration (iii)) are analogous to a double solid twistor spaces on $3\mathbb{CP}^2$ 
of a degenerate form found by Kreussler-Kurke \cite[p.\,50, Case (b)]{KK92}.
Here we mention that there is {\em one more}\, another type of twistor spaces
on $4\mathbb{CP}^2$ of a degenerate form which do not admit $\mathbb C^*$-action.
We will study these 2 kinds of double solid twistor 
spaces  on $4\mathbb{CP}^2$ in a separate paper.
On the other hand, although Moishezon,
it turns out that  the remaining case (iv) does not have structure of double solids,
because the anticanonical map of the twistor spaces
becomes birational.
So they are rather similar to the twistor spaces of Joyce metrics on $4\mathbb{CP}^2$ of non-LeBrun type
\cite{Hon_JAG}.
But contrary to the Joyce's case, explicit 
realization of the anticanonical model seems difficult.

\section{Appendix: Inverting the contraction map $Z_3\to Z_4$
by a blowup}\label{ss:blowup}
We recall from Section \ref{ss:modify} that the singular variety $Z_4$ is obtained from the twistor space $Z$ via
non-singular spaces  $Z_1,Z_2$ and $Z_3$,
and the transformations therein are standard until getting $Z_3$.
On the other hand, the map $\mu_4:Z_3\to Z_4$ contracts the reducible connected divisor $E_1\cup E_1$
to a reducible connected curve $l_4\cup \ol l_4$,
and there we used Fujiki's contraction theorem.
One would wish to find such a birational morphism through the usual procedure of blowups with non-singular center.
In this subsection, we explicitly see that an embedded blowup at  and 
a small resolution provide the desired situation and point out
that the process is in a sense a singular version of the Hironaka's construction of non-projective Moishezon 3-folds.

For this we consider the double covering of $Y$ with branch $B$ which is the intersection with 
the quartic hypersurface \eqref{scrB}.
Recall that $l=\{z_0=z_1=z_2=0\}$, and $B\cap l$ consists of  two points
$\{Q=0\}\cap l$.
As before let $p_2$ be any one of the two points and we work in a neighborhood of $p_2$ as the situation around $\ol p_2$ can see by just taking the image under the real structure.
As in the proof of Proposition \ref{prop:A3},  putting $x=z_0/z_4, y= z_1/z_4, z= z_2/z_4$ and $u=Q$,
we can use $(x,y,z,u)$ as coordinates in a neighborhood of $p_2$ in $ \mathbb{CP}^4$, and
we may suppose that the hyperquartic \eqref{scrB} is defined by  $x=u^2$,
while $Y$ is defined by $x^2=yz$.
Next for studying the structure of the double covering, we introduce another coordinate $w$ over the neighborhood of $p_2$, 
so that the double cover is defined by 
\begin{align}\label{dc10}
x^2=yz, \quad w^2= u^2-x \hsp{\text{in}}\hsp\mathbb C^5 \hsp{\text{with coordinates}}
\hsp (x,y,z,u,w).
\end{align}
(This is the equation of $Z_4$ around the points $l_4\cap \ol l_4$ we promised in Section \ref{ss:modify}.)
Let $W$ be this double covering and $\varpi:W\to Y$  the projection. (Of course this is valid only in a neighborhood of $p_2$.) The singular locus of $W$ is 
\begin{align}\label{2lines}
\varpi^{-1}(l)= \{x=y=z=u-w=0\}\cup \{ x=y=z=u+w=0\},
\end{align}
which is a union of 2 lines, and $W$ has A$_1$-singularities
along these lines minus the origin.
We note that substituting $x=u^2-w^2$ to $x^2=yz$, we obtain that $W$ contains the following distinguished 4 surfaces
\begin{gather}\label{4planes1}
\{ x= y= w-u = 0\},\hsp\{ x= z= w-u = 0\},\\
\label{4planes2}
 \{ x= y= w+u = 0\},\hsp \{ x= z= w+u = 0\}.
\end{gather} 
Then obviously we have 
 $W\cap\{x=y=z=0\}=\varpi^{-1}(l)$.
So if we let $\tilde{\mathbb C}^5\to \mathbb C^5$ to be the blowup at the plane $\{x=y=z=0\}$
and $\tilde W$  to mean the strict transform of $W$,
then $\tilde W\to W$ is an embedded blowup at Sing\,$W$.
Then by concrete computations using coordinates it is not difficult to see that 
the exceptional locus of $\tilde W\to W$ consists of 2 irreducible divisors which are over the 2 lines
\eqref{2lines} respectively,
that the inverse image of the origin is a non-singular rational curve,
and that the singularities of $\tilde W$ consists of 2 points lying on this rational curve, both of which are 
ordinary double points.
Further, by the effect of the blowup, the pair of planes \eqref{4planes1},  both of which contain the same line $\{x=y=z=w-u=0\}$
are separated by one of the exceptional divisors, and the same for another pair \eqref{4planes2}.
This way we get the situation of Figure \ref{fig:trans}, (e).
Then an appropriate small resolution (which is obvious from the figure) 
gives the desired space $Z_3$.

As above the center of the blowup is a reducible curve 
whose fundamental group is $\mathbb Z$.
Thus, together with an inspection of the choice of the  small resolution 
displayed in Figure \eqref{fig:trans}, it would be possible to say that
the transformation from $Z_4$ to $Z_3$ is a singular version
of Hironaka's well-known example of  non-projective Moishezon 3-folds
\cite{Hi60}
in the sense that the center of the blowup in the present situation is a
singular locus of the 3-fold.

\end{document}